\documentclass{amsart}
\usepackage[utf8]{inputenc}
\usepackage[OT2,T1]{fontenc}
\DeclareSymbolFont{cyrletters}{OT2}{wncyr}{m}{n}
\DeclareMathSymbol{\Sha}{\mathalpha}{cyrletters}{"58}

\usepackage{amsmath}
\usepackage{amssymb}
\usepackage{tikz-cd}

\usepackage{amsthm,thmtools}
\usepackage[hidelinks]{hyperref}
\usepackage{cleveref}

\crefname{section}{\S\!\!}{\S\S\!\!}

\theoremstyle{plain}
\newtheorem{theorem}{Theorem}[section]
\newtheorem{lemma}[theorem]{Lemma}
\newtheorem{proposition}[theorem]{Proposition}
\newtheorem{corollary}[theorem]{Corollary}

\theoremstyle{definition}
\newtheorem{definition}[theorem]{Definition}
\newtheorem{example}[theorem]{Example}

\newtheorem{question}{Question}

\theoremstyle{remark}
\newtheorem{remark}[theorem]{Remark}

\title{Arinkin's sheaf and Hilbert schemes}
\date\today
\author{David Zhiyuan Bai}
\address{Yale University}
\email{david.bai@yale.edu}

\begin{document}
\begin{abstract}
    Arinkin's sheaf gives a relative Fourier-Mukai autoequivalence of a general Beauville-Mukai system, extending the usual Poincar\'e bundle over the smooth locus.
    We provide an alternative construction of it as a deformation of the Fourier-Mukai autoequivalence on the Hilbert scheme of an elliptic K3 surface obtained by lifting the classical Poincar\'e kernel for an elliptic surface.

    As applications, we show that the $P=C$ property for a general Beauville-Mukai system is equivalent to that for the Hilbert scheme fibration, thus providing new proofs for both.
    We also construct an explicit deformation that connects Markman's construction of a projectively hyperholomorphic bundle on products of hyper-K\"ahler manifolds of K3${}^{[g]}$-type to that of Hartlieb and Shah.
    This answers a question of Hartlieb and Shah.
\end{abstract}
\maketitle
\tableofcontents
\section{Introduction}
\subsection{Duality}
Let $B$ be a smooth projective base.
Suppose $M\to B$ is an abelian scheme.
The classical Mukai duality supplies a dual abelian scheme $M^\vee\to B$ along with a Poincar\'e line bundle
\[\mathcal L\in\operatorname{Pic}(M\times_B M^\vee)\]
such that the induced Fourier-Mukai transform $\Phi_{\mathcal L}:D^b\operatorname{Coh}(M)\to D^b\operatorname{Coh}(M^\vee)$ is an equivalence.

There are, however, many situations where one has a fibration $M\to B$ that is only generically abelian, i.e.~the generic fiber is an abelian variety, but singular fibers may exist.
Over the smooth locus, we still have a dual with a Poincar\'e bundle (possibly twisted, if a section is not present).
The question then becomes:
\begin{question}
    In what situations does one have a dual $M^\vee\to B$ over the entire base so that the Poincar\'e bundle extends to an object over $M\times_B M^\vee$ inducing a derived equivalence?
\end{question}
\subsection{Arinkin versus BKR}
We are particularly interested in the situation where $M$ is an irreducible holomorphic symplectic variety.
Any fibration of $M$ is Lagrangian \cite{mat_fib}, hence generically abelian.

There are two vastly different ways to construct such a Poincar\'e kernel in this case, both available in examples where $M$ is of K3$^{[g]}$-type.

The first method is essentially due to Arinkin \cite{arinkin}.
It applies to Beauville-Mukai systems for integral curve classes:
Suppose $S$ is a K3 surface and $H$ an integral curve class on $S$ with $H^2=2g-2$.
We can form the moduli space $M_S(0,H,1-g)$ of sheaves on $S$ with Mukai vector $(0,H,1-g)$, which admits a fibration over $|H|\cong\mathbb P^g$.

In this case, the fibers are degree-zero compactified Jacobians of curves of class $H$.
The seminal work of Arinkin \cite{arinkin} constructed a maximal Cohen-Macaulay sheaf
\[\mathcal G_{\rm Arinkin}\in\operatorname{Coh}(M_S(0,H,1-g)\times_{|H|}M_S(0,H,1-g))\]
from the universal sheaf for the compactified Jacobian.

On the other hand, there is a second method that applies to the following rather different situation:
Suppose $S_0\to\mathbb P^1$ is an elliptic K3 surface with integral fibers and a section.
Then there is a fibration $S_0^{[g]}\to S_0^{(g)}\to(\mathbb P^1)^{(g)}=\mathbb P^g$ by composing the symmetric power with the Hilbert-Chow resolution.
Construction of a Poincar\'e kernel for this fibration is currently inaccessible from the compactified Jacobian perspective, as it has non-reduced fibers.

Associated to the elliptic surface $S_0\to\mathbb P^1$ is a Fourier-Mukai autoequivalence $\Phi_{\mathcal P}$ of $D^b\operatorname{Coh}(S_0)$, which can be upgraded to an autoequivalence of $D^b\operatorname{Coh}(S_0^{[g]})$ via the Bridgeland-King-Reid correspondence \cite{bkr}.
We call the result the BKR congruence.
See \cref{s:actionBKR} for the precise construction.
By \cite{bai25}, its kernel is a maximal Cohen-Macaulay sheaf
\[\mathcal P^{[g]}\in\operatorname{Coh}(S_0^{[g]}\times_{\mathbb P^g}S_0^{[g]}).\]
Unlike the compactified Jacobian approach, this utilizes the geometry of Hilbert schemes in lieu of understanding the fibers of the fibration.

In this article, we develop a relationship between these two constructions.
We utilize a degeneration of Beauville-Mukai systems for general polarized K3 surfaces to this Hilbert scheme fibration due to Yoshioka \cite{yoshioka}.
Under this degeneration, we show:
\begin{theorem}
    The BKR congruence $\mathcal P^{[g]}$ is a specialization of $\mathcal G_{\rm Arinkin}^\vee$.
\end{theorem}
This is obtained by deforming $\mathcal P^{[g]}$ away using the semiregularity theorem of Buchweitz and Flenner \cite{semireg}.
The kernel of a derived equivalence between Calabi-Yau varieties is always semiregular \cite{rienks}, so the main input is to show that its Chern character remains Hodge along the deformation.

We show this using the theory of extended Mukai lattice developed by Taelman \cite{taelman} and Beckmann \cite{beckmann}.
This allows us to reduce the problem of showing $\operatorname{ch}(\mathcal P^{[g]})$ remains Hodge to showing that the induced action of $\Phi_{\mathcal P^{[g]}}$ on the extended Mukai lattice $\tilde{H}(S_0^{[g]};\mathbb Q)$ remains Hodge, which is much simpler.
This is done in \cref{s:yoshioka}.

We then show that the deformation obtained is locally unique, and is a maximal Cohen-Macaulay sheaf supported on $M_S(0,H,1-g)\times_{|H|}M_S(0,H,1-g)$ inducing a derived autoequivalence.
As computed in \cref{s:comparison}, it is precisely $\mathcal G^\vee_{\rm Arinkin}$.

Note that we do not, a priori, need to use the existence of Arinkin's sheaf to obtain this locally unique deformation as a maximal Cohen-Macaulay Fourier-Mukai kernel.
It therefore gives another proof of:
\begin{theorem}[\cite{arinkin}]
    For a general polarized K3 surface $(S,H)$ of genus $g$, there exists a maximal Cohen-Macaulay extension of the usual Poincar\'e bundle over the smooth locus to
    \[M_S(0,H,1-g)\times_{|H|}M_S(0,H,1-g),\]
    and it induces a Fourier-Mukai autoequivalence.
\end{theorem}
In other words, it gives a different construction of Arinkin's sheaf.
On the level of curves, this amounts to degenerating an integral curve of arithmetic genus $g$ into a (reducible, and possibly non-reduced) centipede-shaped curve of the form
\[s+F_1+\cdots+F_g\]
where $s,F_j$ are the section and elliptic fibers of $S_0\to\mathbb P^1$.
The Poincar\'e kernel for this degenerate situation is supplied by the Poincar\'e kernel for the elliptic surface, and Arinkin's sheaf is constructed by deforming it away.

In the future, we hope to adapt similar methods to construct Poincar\'e kernels in other situations, where Arinkin's construction may be unavailable.

\subsection{An application: the \texorpdfstring{$P=C$}{P=C} property}
For a fibration $h:M\to B$ with $M$ smooth projective, one can define the perverse filtration on the cohomology of $M$ by setting
\[P_kH^l(M;\mathbb Q)=\operatorname{Im}(H^l(B,{}^{\mathfrak p}\tau_{\le k+\dim B}\mathbf Rh_\ast\underline{\mathbb Q})\to H^l(B,\mathbf Rh_\ast\underline{\mathbb Q})=H^l(M;\mathbb Q)),\]
where ${}^{\mathfrak p}\tau$ are perverse truncation functors.

The $P=C$ property, as studied in depth in the work of Maulik, Shen, and Yin \cite{msy}, roughly says that the perverse filtration for $M\to B$, which is entirely topological, can be recovered algebraically through the Chern classes of a certain object $\mathcal G$ in the derived category of coherent sheaves.
See \cref{s:comparison} for a precise definition.

As a consequence of our degeneration, we show (\Cref{prop:cohP=C}) that the $P=C$ property for a general Beauville-Mukai system $M_S(0,H,1-g)\to |H|$ with $\mathcal G=\mathcal G_{\rm Arinkin}$ is equivalent to that for the Hilbert scheme fibration $S_0^{[g]}\to\mathbb P^g$ with $\mathcal G=\mathcal P^{[g]}$.
This result can also be categorified (\Cref{prop:P=Ccat}).

By \cite{msy}, $P=C$ does hold for general Beauville-Mukai systems, so this establishes $P=C$ for the Hilbert scheme fibration.

Conversely, upcoming work of Ghosh \cite{soumik} proves $P=C$ for the Hilbert scheme fibration using Nakajima operators.
The theorem therefore provides a different proof of $P=C$ for general Beauville-Mukai systems.

\subsection{Markman versus Hartlieb-Shah}
The degeneration also helps compare two different constructions of projectively hyperholomorphic bundles on products of hyper-K\"ahler manifolds of K3$^{[g]}$-type.

In \cite{markmanhyperhol}, Markman started with a general polarized K3 surface $(S,H)$ and considered an isotropic Mukai vector of the form
\[(n,H,s),n\ge 2,\gcd(n,s)=1.\]
The associated moduli space is then another K3 surface which is derived equivalent to $S$, with kernel given by the universal bundle.
The Bridgeland-King-Reid equivalence then allows one to lift it to a Fourier-Mukai equivalence between the relevant Hilbert schemes.

Markman showed that the kernel of this lift is a projectively hyperholomorphic vector bundle.
The existence of such a bundle has various applications, such as the algebraicity of Hodge isometries in the same work \cite{markmanhyperhol}, as well as the $D$-equivalence conjecture \cite{msyz}.

Recently, Hartlieb and Shah \cite{hartliebshah} provided a very different construction of a projectively hyperholomorphic bundle.
They instead considered a hyper-K\"ahler variety $Y$ of K3${}^{[g]}$-type with two transverse Lagrangian fibrations
\[\begin{tikzcd}
    X\drar&&Y\dlar\drar&&Z\dlar\\
    &\mathbb P^g&&\mathbb P^g&
\end{tikzcd},\]
each being a degenerate twistor line deformation of very general Beauville-Mukai systems $X\to\mathbb P^g,Z\to\mathbb P^g$.
Arinkin's sheaf then provides twisted equivalences
\[D^b(Y)\cong D^b(X,\alpha_X),D^b(Y)\cong D^b(Z,\alpha_Z).\]

They showed that the composition of these two equivalences provides a projectively hyperholomorphic (twisted) bundle on $X\times Z$.
%Let us note that both constructions produce rigid (twisted) bundles (\Cref{lem:FMkernelsarerigid}).

In \cite[Question 1.6]{hartliebshah}, Hartlieb and Shah asked whether their construction deforms to an instance of Markman's under diagonal twistor line deformations.
In \cref{s:hyperhol}, we give an affirmative answer:
\begin{theorem}
    There exist diagonal twistor line deformations that connect the construction of Hartlieb and Shah to that of Markman.
\end{theorem}
In fact, we show the following stronger statement:
Let $S_0$ be a K3 surface with two elliptic fibrations $\pi_i$, both with integral fibers and a section.
And let $\mathcal P_i$ be the kernel for the relative autoequivalence associated to $\pi_i$.
Then the composition
\[\mathcal P_1^{[g]}\circ(\mathcal P_2^{[g]})^{-1}[-g]\]
is a rigid, projectively hyperholomorphic bundle, and is simultaneously a specialization of families of Hartlieb-Shah's constructions (\Cref{cor:hartliebshah_specialize}) and families of Markman's constructions (\Cref{prop:markmanspecialize}).

Along the way, we also provide a different method to construct the space $Y$ in Hartlieb-Shah's argument (\Cref{lem:reconstr_hartliebshah}), as well as a different proof of theorem of the square for $\mathcal P^{[g]}$ as in \cite[Theorem 0.3]{bai25} for K3 surfaces (\Cref{rem:altproof_tots}).
\subsection*{Acknowledgement}
The author is very grateful to Junliang Shen, who suggested to the author that such a specialization may be possible.
The author also benefitted from discussions with David Fang and Soumik Ghosh.

This work is partially supported by NSF grant DMS-2301474.
\subsection*{AI disclosure}
An AI tool was used to proofread the manuscript.
It pointed out a gap in a previous version of \Cref{aicorrectedthis}, where the author prematurely claimed that $f=\operatorname{id}$.
\section{Action of the BKR congruence}\label{s:actionBKR}
\subsection{The BKR congruence}
Let $S$ be a K3 surface.
The symmetric group $\operatorname{Sym}_g$ on $g$ letters acts on the product $S^g$ by permuting the coordinates.
Let $\tilde{S}_g\hookrightarrow S^{[g]}\times S^g$ be the $\operatorname{Sym}_g$-isospectral Hilbert scheme.
It comes with natural projections
\[\begin{tikzcd}
    S^{[g]}&{[\tilde{S}_g/\operatorname{Sym}_g]}\rar{p}\lar[swap]{q}&{[S^g/\operatorname{Sym}_g]}
\end{tikzcd}.\]
Note that $p$ is birational while $q$ is a $\operatorname{Sym}_g$-quotient.
\begin{theorem}[\cite{bkr}]
    The functor
    \[\operatorname{BKR}=\mathbf Rp_\ast q^\ast:D^b\operatorname{Coh}(S^{[g]})\to D^b\operatorname{Coh}([S^g/\operatorname{Sym}_g])\]
    is an equivalence of categories.
\end{theorem}
If one considers the kernel the other way around, one obtains Krug's functor \cite{krug}:
\[\operatorname{BKR}^\top=q_\ast\mathbf L p^\ast:D^b\operatorname{Coh}([S^g/\operatorname{Sym}_g])\to D^b\operatorname{Coh}(S^{[g]}),\]
which is also an equivalence.

These functors allow one to lift an equivalence between K3 surfaces to their Hilbert schemes.
For simplicity, we only consider autoequivalences.
Let $\Phi:D^b\operatorname{Coh}(S)\to D^b\operatorname{Coh}(S)$ be a Fourier-Mukai equivalence with kernel $\mathcal K$.
Then we can view $\mathcal K^{\boxtimes g}$ as an object in $D^b\operatorname{Coh}([(S^g\times S^g)/\operatorname{Sym}_g^{\rm diag}])$, where $\operatorname{Sym}_g^{\rm diag}$ is the diagonal copy of $\operatorname{Sym}_g$ in $\operatorname{Sym}_g\times \operatorname{Sym}_g$.
It induces an autoequivalence:
\[\Phi^{\boxtimes g}=\Phi_{\mathcal K^{\boxtimes g}}:D^b\operatorname{Coh}([S^g/\operatorname{Sym}_g])\to D^b\operatorname{Coh}([S^g/\operatorname{Sym}_g]).\]

Through this and BKR, one can then lift $\Phi$ to a Fourier-Mukai autoequivalence for $S^{[g]}$, which we shall refer to as the BKR congruence:
\[\Phi^{[g]}=\operatorname{BKR}^\top\circ\Phi^{\boxtimes g}\circ\operatorname{BKR}:D^b(S^{[g]})\to D^b(S^{[g]}).\]
Let us also define the BKR conjugate:
\[\Phi^{[g]+}=\operatorname{BKR}^\top\circ\Phi^{\boxtimes g}\circ(\operatorname{BKR}^\top)^{-1}=\Phi^{[g]}\circ\operatorname{BKR}^{-1}\circ(\operatorname{BKR}^\top)^{-1}.\]

We write $\mathcal K^{[g]}$ and $\mathcal K^{[g]+}$ for the kernels of these equivalences, respectively.

Let $\chi:D^b_{\operatorname{Sym}_g}(S^g)\to D^b_{\operatorname{Sym}_g}(S^g)$ be the sign involution, and let $\delta$ be half the exceptional divisor of the Hilbert-Chow resolution $S^{[g]}\to S^{(g)}$.
Write
\[\chi^{[g]}=\operatorname{BKR}^\top\circ\chi\circ(\operatorname{BKR}^\top)^{-1}.\]
The following computation is due to Markman \cite{markmanBKR}.
\begin{lemma}[{\cite[(11.5)]{markmanBKR}}]\label{lemma:congruence_conjugate_factorize}
    We have
    \[\operatorname{BKR}^\top\circ \operatorname{BKR}=\chi^{[g]}\circ(-\otimes\mathcal O(-\delta)).\]
    In particular, $\Phi^{[g]}=\Phi^{[g]+}\circ\chi^{[g]}\circ (-\otimes\mathcal O(-\delta))$.
\end{lemma}

Now let $\pi:S_0\to\mathbb P^1$ be an elliptic K3 surface with integral fibers and a section $s:\mathbb P^1\to S_0$.
Then one has a relative Poincar\'e duality
\[\mathcal F:D^b\operatorname{Coh}(S_0)\to D^b\operatorname{Coh}(S_0)\]
given by the kernel
\[\mathcal P=\mathcal I_{\Delta/(S_0\times_{\mathbb P^1}S_0)}\otimes\operatorname{pr}_1^\ast\mathcal O(s)\otimes\operatorname{pr}_2^\ast\mathcal O(s)\otimes\operatorname{pr}_{\mathbb P^1}^\ast\mathcal O(2)\in \operatorname{Coh}(S_0\times_{\mathbb P^1}S_0).\]
This normalization is chosen so that $\mathcal F(\mathcal O_s)=\mathcal O_{S_0}$.

We may therefore form the BKR congruence for this equivalence.
Since the kernel $\mathcal P$ is a maximal Cohen-Macaulay sheaf on $S_0\times_{\mathbb P^1}S_0$, the kernel of $\mathcal F^{[g]}$ inherits some of these properties.

To state them, let us consider the fibration $\pi^{[g]}:S_0^{[g]}\to\mathbb P^g$ given by composing the Hilbert-Chow resolution $S_0^{[g]}\to S_0^{(g)}$ with the natural map $S_0^{(g)}\to (\mathbb P^1)^{(g)}=\mathbb P^g$.
There is a section $s^{[g]}:\mathbb P^g\to S_0^{[g]}$ from $s:\mathbb P^1\to S_0$.
\begin{theorem}[{\cite[Corollary 3.12, Proposition 4.1]{bai25}}]
    The kernel $\mathcal P^{[g]}$ of $\mathcal F^{[g]}$ is a maximal Cohen-Macaulay sheaf supported on $S_0^{[g]}\times_{\mathbb P^g}S_0^{[g]}$ and is invariant under swapping the factors.
    Moreover, $\mathcal P^{[g]}|_{s^{[g]}\times_{\mathbb P^g}S_0^{[g]}}=\mathcal O_{s^{[g]}\times_{\mathbb P^g}S_0^{[g]}}$.
\end{theorem}

For later convenience, let us compute the action of $\mathcal F$ on the extended Mukai lattice of $S_0$.
Write $s,f$ for the fundamental classes of the section $s$ and of a fiber, respectively.
Then $\tilde{H}(S_0;\mathbb Q)=\mathbb Q\alpha\oplus H^2(S_0;\mathbb Q)\oplus\mathbb Q\beta$, where $\beta$ is the class of a point, and $\alpha$ is the class of $S_0$.
\begin{proposition}\label{prop:elliptic_surface_calc}
    On $\tilde{H}(S_0;\mathbb Q)$, the cohomological realization $\mathcal F^{\tilde{H}}$ of $\mathcal F$ acts by
    \[\mathcal F^{\tilde{H}}(\beta)=f,\mathcal F^{\tilde{H}}(f)=-\beta,\mathcal F^{\tilde{H}}(s)=\alpha-f+\beta,\mathcal F^{\tilde{H}}(\alpha)=-s-f+\beta\]
    and $\mathcal F^{\tilde{H}}(c)=-c$ for any $c\in \langle s,f\rangle_{H^2}^\perp$.
\end{proposition}
Note, in particular, that $\det \mathcal F^{\tilde{H}}=1$.
\begin{proof}
    Let $f$ also denote a smooth fiber of $\pi$ and let $p$ be the intersection of $f$ and $s$.
    Then
    \[\mathcal F(k_p)=\mathcal O_f,\mathcal F(\mathcal O_f)=k_p[-1],\mathcal F(\mathcal O_s)=\mathcal O_{S_0},\mathcal F(\mathcal O_{S_0})=\mathcal O_s(-2)[-1].\]
    Passing these to Mukai vectors shows the calculations for the $\beta,f,s,\alpha$ classes.

    Now let $c\in\langle s,f\rangle_{H^2}^\perp$.
    Factorize
    \[\mathcal P=(-\otimes\mathcal O(s+f))\circ\mathcal I_{\Delta/(S_0\times_{\mathbb P^1}S_0)}\circ (-\otimes\mathcal O(s+f)).\]
    Since $c\cdot s=c\cdot f=0$, $c$ is a fixed point for the action of $(-\otimes\mathcal O(s+f))$.
    So it remains to show that the action induced by $\mathcal I_{\Delta/(S_0\times_{\mathbb P^1}S_0)}$ sends $c$ to $-c$.

    In the K-theory of $S_0\times S_0$, we have
    \[[\mathcal I_{\Delta/(S_0\times_{\mathbb P^1}S_0)}]=[S_0\times_{\mathbb P^1}S_0]-[\Delta]=[S_0\times S_0]-[\mathcal O(-f)\boxtimes\mathcal O(-f)]-[\Delta].\]
    But for a line bundle $L\in\operatorname{Pic}(S_0)$, the action of $[L\boxtimes L]$ sends $c$ to $(c\cdot c_1(L))v(L)$ which is zero for $L=\mathcal O,\mathcal O(-f)$.
    Hence $[\mathcal I_{\Delta/(S_0\times_{\mathbb P^1}S_0)}]$ sends $c$ to $-c$.
\end{proof}

\subsection{The extended Mukai lattice}
Let us start by reviewing the theory of the LLV algebra acting on the cohomology of a compact hyper-K\"ahler variety.

For a compact hyper-K\"ahler variety $X$ of dimension $2d$, we say a class $\lambda\in H^2(X;\mathbb Q)$ enjoys the hard Lefschetz property if the operator $e_\lambda=-\smile \lambda$ has the property that
\[e_\lambda^i:H^{2d-i}(X;\mathbb Q)\to H^{2d+i}(X;\mathbb Q)\]
is an isomorphism for all $i>0$.
The Jacobson-Morosov theorem allows us to complete any such operator to an $\mathfrak{sl}_2$-triple $(e_\lambda,f,h)$ where $h$ is the cohomology grading $h|_{H^k}=(k-2d)\operatorname{id}$.
\begin{definition}
    The Looijenga-Lunts-Verbitsky (LLV) algebra $\mathfrak g(X)$ of $X$ is the sub-Lie algebra of $\operatorname{End}(H^\ast(X;\mathbb Q))$ generated by all such $\mathfrak{sl}_2$-triples.
\end{definition}
Note that the condition for a class to have the hard Lefschetz property is Zariski open, so these classes span $H^2(X;\mathbb Q)$.
In particular, $e_\lambda\in \mathfrak g(X)$ for all $\lambda\in H^2(X;\mathbb Q)$.

The LLV algebra can be described by introducing the extended Mukai lattice:
\begin{definition}
    The extended Mukai lattice of $X$ is the lattice
    \[(\mathbb Q\alpha\oplus H^2(X;\mathbb Q)\oplus\mathbb Q\beta,b)\]
    where $b|_{H^2(X;\mathbb Q)}$ is the Beauville-Bogomolov-Fujiki (BBF) form, and
    \[b(\alpha,\alpha)=0,b(\beta,\beta)=0,b(\alpha,\beta)=-1,H^2(X;\mathbb Q)\perp \mathbb Q\alpha\oplus\mathbb Q\beta.\]
    It is given the structure of a weight $2$ Hodge structure by taking the usual weight $2$ Hodge structure on $H^2(X;\mathbb Q)$ and setting $\alpha,\beta$ at type $(1,1)$.
\end{definition}
We will sometimes write $x\cdot y$ in place of $b(x,y)$.

For the definition and basic properties of the BBF form, we refer the reader to \cite{verbitsky}.
Note that any $\lambda\in H^2(X;\mathbb Q)$ defines an operator $e_\lambda:\tilde{H}(X;\mathbb Q)\to \tilde{H}(X;\mathbb Q)$ sending
\[\alpha\mapsto\lambda,\lambda'\mapsto (\lambda\cdot\lambda')\beta,\beta\mapsto 0\]
for any $\lambda'\in H^2(X;\mathbb Q)$.
\begin{theorem}[{\cite[Theorem 11.1]{verbitsky}}]
    There is an isomorphism $\mathfrak g(X)\cong\mathfrak{so}(\tilde{H}(X;\mathbb Q))$ that sends $e_\lambda$ to $e_\lambda$ for any class $\lambda\in H^2(X;\mathbb Q)$.
\end{theorem}
This isomorphism may be understood through the Verbitsky component of cohomology.
\begin{definition}
    The Verbitsky component $\operatorname{SH}(X;\mathbb Q)$ of $H^\ast(X;\mathbb Q)$ is the subalgebra generated by $H^2(X;\mathbb Q)$.
\end{definition}
It is a sub-$\mathfrak g(X)$-module of $H^\ast(X;\mathbb Q)$ (see e.g.~\cite[Proof of Theorem 2.15]{GKLR}).

Consider the natural map
\[\tilde\Psi:\operatorname{Sym}^\bullet H^2(X;\mathbb Q)\to\operatorname{Sym}^d\tilde{H}(X;\mathbb Q),\lambda_1\cdots\lambda_n\mapsto e_{\lambda_1}\cdots e_{\lambda_n}\frac{\alpha^d}{d!}\]
where the action of $e_\lambda$ on $\tilde{H}(X;\mathbb Q)$ is extended to $\operatorname{Sym}^d\tilde{H}(X;\mathbb Q)$ by setting
\[e_\lambda(\gamma_1\cdots\gamma_d)=\sum_{j}\gamma_1\cdots e_\lambda(\gamma_j)\cdots\gamma_d.\]
\begin{theorem}[{\cite[Proposition 3.5, Lemma 3.7]{taelman}}]\label{thm:laplacian}
    The map $\tilde\Psi$ factors through a $\mathfrak g(X)$-equivariant injection
    \[\Psi:\operatorname{SH}(X;\mathbb Q)\to\operatorname{Sym}^d\tilde{H}(X;\mathbb Q)\]
    whose cokernel is the Laplacian $\Delta:\operatorname{Sym}^d\tilde{H}(X;\mathbb Q)\to \operatorname{Sym}^{d-2}\tilde{H}(X;\mathbb Q)$ which sends
    \[\gamma_1\cdots\gamma_d\mapsto\sum_{i<j}b(\gamma_i,\gamma_j)\gamma_1\cdots\hat{\gamma}_i\cdots\hat{\gamma}_j\cdots\gamma_d.\]

    In particular, $\operatorname{SH}(X;\mathbb Q)$ is a faithful $\mathfrak g(X)$-module.
\end{theorem}
% For any $\mathfrak g(X)$-equivariant map $\phi:\tilde{H}(X;\mathbb Q)\to \tilde{H}(X;\mathbb Q)$, we then get a map
% \[S_{[d]}(\phi):\operatorname{SH}(X;\mathbb Q)\to\operatorname{SH}(X;\mathbb Q)\]
% from the identification of $\operatorname{SH}(X;\mathbb Q)$ as the kernel of the Laplacian.
\begin{theorem}[{\cite[Proposition 4.2]{taelman}}]
    Suppose $f:H^\ast(X;\mathbb Q)\to H^\ast (X;\mathbb Q)$ is an isomorphism preserving $\operatorname{SH}(X;\mathbb Q)$ such that
    \[f^{-1}\mathfrak g(X) f=\mathfrak g(X),\]
    then:
    \begin{enumerate}
        \item there exists a scalar $\lambda$ and a similitude $\phi:\tilde{H}(X;\mathbb Q)\to \tilde{H}(X;\mathbb Q)$ such that the diagram
        \[
        \begin{tikzcd}
            \operatorname{SH}(X;\mathbb Q)\rar{f}\dar{\Psi}&\operatorname{SH}(X;\mathbb Q)\dar{\Psi}\\
            \operatorname{Sym}^d\tilde{H}(X;\mathbb Q)\rar{\lambda \operatorname{Sym}^d\phi}&\operatorname{Sym}^d\tilde{H}(X;\mathbb Q)
        \end{tikzcd}
        \]
        commutes, moreover,
        \item the pair $(\lambda,\phi)$ is essentially unique, in the sense a different choice $(\lambda',\phi')$ necessarily has the form $(\lambda\mu^{-d},\mu\phi)$ for some nonzero scalar $\mu$.
    \end{enumerate}
\end{theorem}
\begin{corollary}\label{cor:commutewllv}
    For any $g\in\mathfrak g(X)=\mathfrak{so}(\tilde{H}(X;\mathbb Q))$ and any choice of $\phi$ in the preceding theorem, we have
    \[f^{-1}gf=\phi^{-1}g\phi.\]
\end{corollary}
\begin{proof}
    Since $\operatorname{SH}(X;\mathbb Q)$ is a faithful $\mathfrak g(X)$-representation, it suffices to show that the two sides have the same action on $\operatorname{SH}(X;\mathbb Q)$.
    But $f$ is the restriction of $\lambda\operatorname{Sym}^d\phi$ by the preceding theorem, so this is clear.
\end{proof}
We write $f^{\tilde{H}}$ to denote such a $\phi$ up to a nonzero scalar.
When $f$ is induced by a derived equivalence, we can sometimes make a distinguished choice for $f^{\tilde{H}}$:
\begin{theorem}[{\cite[\S4]{taelman}}]
    Suppose $\Phi:D^b(X)\to D^b(X)$ is a derived equivalence and write $f=\Phi^H:H^\ast(X;\mathbb Q)\to H^\ast(X;\mathbb Q)$.
    Then:
    \begin{enumerate}
        \item $f^{-1}\mathfrak g(X)f=\mathfrak g(X)$,
        \item $f$ preserves $\operatorname{SH}(X;\mathbb Q)$,
        \item the scalar $\lambda$ can be chosen to be $1$ if $d$ is odd, and
        \item the scalar $\lambda$ can be uniquely chosen to be one of $\{\pm 1\}$ if $d$ is even and $b_2(X)$ is odd.
    \end{enumerate}
    Moreover, $\phi$ respects the weight $2$ Hodge structure on $\tilde{H}(X;\mathbb Q)$.
\end{theorem}
Note that $b_2(X)=23$ is odd if $X$ is of K3$^{[d]}$-type for $d>1$.
We write $\Phi^{\tilde{H}}$ for the choice of $\phi$ corresponding to the choice of $\lambda$ in either the situation of (3) or (4).
In situation (4), we require $\lambda=\det\Phi^{\tilde{H}}$.

The formation of $\Phi^{\tilde{H}}$ is then functorial in $\Phi$.
\subsection{Computation of the induced action}
Let us now compute the action of $(\mathcal F^{[g]})^{\tilde{H}}$.
The main ingredient is Beckmann's calculation of the extended Mukai action of various BKR conjugates.
We will only be interested in hyper-K\"ahler varieties $X$ of K3$^{[g]}$-type.

For a K3 surface $S$, let us make the identifications
\[H^2(S^{[g]},\mathbb Q)\cong H^2(S,\mathbb Q)\oplus\mathbb Q\delta,\tilde{H}(S^{[g]},\mathbb Q)\cong\tilde{H}(S,\mathbb Q)\oplus\mathbb Q\delta\]
where $\delta$ is half the class of the exceptional divisor.
Note that $b(\delta,\delta)=2-2g$.

Following Beckmann \cite{beckmann}, we define:
\begin{definition}
    For any $\lambda\in H^2(X;\mathbb Q)$, its associated B-field twist
    \[B_\lambda:\tilde{H}(X;\mathbb Q)\to\tilde{H}(X;\mathbb Q)\]
    is defined by
    \[B_\lambda(r\alpha+\mu+s\beta)=r\alpha+\mu+r\lambda+\left(s+b(\lambda,\mu)+r\frac{b(\lambda,\lambda)}{2}\right)\beta\]
    for any $r,s\in\mathbb Q,\mu\in H^2(X;\mathbb Q)$.
\end{definition}
\begin{remark}
    By \cite[Proposition 3.2]{taelman}, $(-\otimes L)^{\tilde{H}}=B_{c_1(L)}$ for any $L\in\operatorname{Pic}(X)$.
\end{remark}
Suppose now that $X=S^{[g]}$ for a K3 surface $S$.
Let us review some calculations done in \cite{beckmann}.
For $\lambda\in\tilde{H}^{1,1}(X;\mathbb Q)$ with $b(\lambda,\lambda)\neq 0$, we write
\[s_\lambda(\mu)=\mu-2\frac{b(\mu,\lambda)}{b(\lambda,\lambda)}\lambda\]
for the reflection across $\lambda$.
\begin{theorem}[{\cite[Proposition 7.1]{beckmann}}]
    Suppose $X=S^{[g]}$ for a K3 surface $S$.
    Then
    \[(\chi^{[g]})^{\tilde{H}}=(-1)^{g+1}s_{\tilde\delta}.\]
    Here, $\tilde\delta=B_{-\delta/2}(\delta)=\delta-(1-g)\beta$.
\end{theorem}
\begin{theorem}[{\cite[Theorem 7.4]{beckmann}}]
    Suppose $\Phi:D^b\operatorname{Coh}(S)\to D^b\operatorname{Coh}(S)$ is a derived equivalence, then
    \[(\Phi^{[g]+})^{\tilde{H}}=(\det\Phi^{\tilde{H}})^{g+1}B_{-\delta/2}\iota(\Phi^{\tilde{H}})B_{\delta/2}.\]
    Here, $\iota(g)(\mu+r\delta)=g(\mu)+r\delta$ for any $\mu\in \tilde{H}(S,\mathbb Q),r\in\mathbb Q$.
\end{theorem}
Let us now take $S=S_0$ and $\Phi=\mathcal F$ the Fourier-Mukai autoequivalence as considered before.
\begin{proposition}\label{prop:actasscalaronperp}
    Let $V=\langle\alpha,s-\delta/2,f,\beta\rangle$.
    Then
    \[(\mathcal F^{[g]})^{\tilde{H}}|_{V^\perp}=(-1)^{g}\operatorname{id}.\]
\end{proposition}
The reason for the choice of $V$ will become apparent later, when we consider a deformation of the Hilbert scheme fibration.
\begin{proof}
    Let us first compute $V^\perp$.
    Since $\langle\alpha,\beta\rangle$ is nondegenerate, $V^\perp\subset\langle s-\delta/2,f\rangle_{H^2}^\perp$.

    But $\langle s-\delta/2,f\rangle$ is also nondegenerate, so $V^\perp$ has dimension $21$.
    Clearly $V^\perp$ contains $\langle s,f\rangle^{\perp}_{H^2(S_0;\mathbb Q)}$ as well as the class $\delta+(1-g)f$.
    They already give $21$ dimensions.
    Therefore
    \[V^\perp=\langle s,f\rangle^{\perp}_{H^2(S_0;\mathbb Q)}\oplus \mathbb Q(\delta+(1-g)f).\]
    We will evaluate their images using the preceding theorems as well as \Cref{prop:elliptic_surface_calc}.

    By \Cref{lemma:congruence_conjugate_factorize} and the preceding theorems, we get
    \begin{align*}
        (\mathcal F^{[g]})^{\tilde{H}}&=(\mathcal F^{[g]+}\circ\chi^{[g]}\circ(-\otimes\mathcal O(-\delta)))^{\tilde{H}}\\
        &=B_{-\delta/2}\iota(\mathcal F^{\tilde{H}})B_{\delta/2}(-1)^{g+1}s_{\tilde\delta}B_{-\delta}.
    \end{align*}

    One readily sees that any $c\in \langle s,f\rangle^{\perp}_{H^2(S_0;\mathbb Q)}$ is a fixed point for $B_{\pm\delta/2},B_{-\delta},s_{\tilde\delta}$, and $\iota(\mathcal F^{\tilde{H}})$ sends it to $-c$.
    So $c$ gets sent to $(-1)^gc$.

    On the other hand, the images of $\delta+(1-g)f$ through various transforms are
    \[\begin{tikzcd}
        \delta+(1-g)f\rar[mapsto]{B_{-\delta}}&\delta+(1-g)f-2(1-g)\beta\dar[mapsto]{s_{\tilde\delta}}\\
        -\delta-(1-g)\beta+(1-g)f\dar[mapsto]{\iota(\mathcal F^{\tilde{H}})}&-\delta+(1-g)f\lar[mapsto]{B_{\delta/2}}\\
        -\delta-(1-g)f-(1-g)\beta\rar[mapsto]{B_{-\delta/2}}&-\delta-(1-g)f
    \end{tikzcd}.\]
    Therefore $\delta+(1-g)f$ gets sent to $(-1)^{g+1}(-\delta-(1-g)f)=(-1)^g(\delta+(1-g)f)$.
\end{proof}
\section{The Yoshioka deformation}\label{s:yoshioka}
\subsection{Yoshioka's isomorphism}
The fibration $S_0^{[g]}\to\mathbb P^g$ admits a different description as a moduli space of one-dimensional sheaves.

For a K3 surface $S$ and $v=(r,c_1,\chi-r)\in\tilde{H}(S;\mathbb Z)$ a primitive Mukai vector, we write $M_S(r,c_1,\chi-r)$ for the moduli space of stable sheaves of Mukai vector $(r,c_1,\chi-r)$ with respect to a generic polarization.
When $r=0$, this moduli space admits a map
\[M_S(0,c_1,\chi)\to |c_1|\]
known as the Beauville-Mukai system.
Here, $|c_1|$ denotes the linear system associated to the (unique) line bundle whose first Chern class is $c_1$.

\begin{theorem}[{\cite[Theorem 3.15]{yoshioka}}]\label{thm:yoshiokaiso}
    There is an isomorphism
    \[S_0^{[g]}\cong M_{S_0}(0,s+gf,1-g)\]
    where the choice of the generic polarization on the right-hand side is $s+kf$ for large $k$.
    This isomorphism is given by sending the ideal $\mathcal I_Z$ of a length $g$ subscheme of $S_0$ to
    \[\mathbf R\underline{\operatorname{Hom}}(\mathcal F(\mathcal I_Z),\mathcal O_{S_0})(-gf).\]

    Moreover, there is an identification $|s+gf|\cong\mathbb P^g$ making the diagram
    \[\begin{tikzcd}
        S_0^{[g]}\rar{\cong}\dar& M_{S_0}(0,s+gf,1-g)\dar\\
        \mathbb P^g=(\mathbb P^1)^{(g)}\rar{\cong}& {|s+gf|}
    \end{tikzcd}\]
    commute.
\end{theorem}
\begin{proof}
    By \cite[Theorem 3.15]{yoshioka}, there is an isomorphism
    \[M_{S_0}(1,(g-2)f,\chi(\mathcal O_{S_0})-1-g)\to M_{S_0}(0,s+gf,g-1),E\mapsto \mathcal F(E)(2f)[1].\]
    On the other hand, $M_{S_0}(1,(g-2)f,\chi(\mathcal O_{S_0})-1-g)$ consists precisely of sheaves of the form $\mathcal I_Z((g-2)f)$ where $Z$ is a length $g$ subscheme of $S_0$.
    We therefore get an isomorphism
    \[S_0^{[g]}\to M_{S_0}(0,s+gf,g-1),\mathcal I_Z\mapsto \mathcal F(\mathcal I_Z)(gf)[1].\]
    Composing this with the isomorphism in \cite[Remark 3.5]{yoshioka} yields the proposed isomorphism.

    As for the identification of the linear system, note simply that $|s+gf|$ consists of effective Cartier divisors of the form $s+F_1+\cdots+F_g$ where each $F_i$ is a fiber of $\pi$.
    So we may identify $|s+gf|\cong(\mathbb P^1)^{(g)}$ by sending $s+F_1+\cdots+F_g$ to $\pi(F_1)+\cdots+\pi(F_g)$.
    Commutativity of the diagram then follows from the calculation in \cite[Eq.~(3.15)]{yoshioka}.
\end{proof}
Suppose now that there is a family of (quasi)-polarized K3 surfaces
\[(\mathcal S,\mathcal H)=\{(S_t,H_t)\}_{t\in\Omega}\to \Omega\]
where $\Omega$ is a small disk, so that $(S_0,H_0)=(S_0,s+gf)$ and $\operatorname{Pic}(S_\eta)=\mathbb ZH_\eta$.
Such a family can be constructed, for example, by putting $(S_0,s+gf)$ in the moduli space of primitively (quasi-)polarized K3 surfaces of genus $g$ and considering a very general deformation.

After possibly shrinking $\Omega$, we may assume that $|H_t|$ consists only of integral curves for all $t\neq 0$.
In particular, the space $M_{S_t}(0,H_t,1-g)$ is just the space of simple sheaves of Mukai vector $(0,H_t,1-g)$, i.e.~any choices of a stability condition are equivalent, for $t\neq 0$.

Let us then form the relative moduli space
\[M_{\mathcal S/\Omega}(0,\mathcal H,1-g)\to \Omega\]
whose fiber over $t\in\Omega$ is $M_{S_t}(0,H_t,1-g)$.
\begin{remark}\label{rem:overKg}
    For $g>2$, the existence of this relative moduli space follows directly from \cite[Proposition 5.1]{yoshioka} since $s+gf$ is ample.
    Let us note that the argument therein works without this assumption (that is, it is also true for $g=2$).
    In fact, a relative moduli space can be constructed without assuming either $\Omega$ is one-dimensional nor $\operatorname{Pic}(S_\eta)=\mathbb ZH_\eta$.
    Indeed, by \cite{kosarew_okonek} one can always form the (possibly non-separated) relative moduli space of simple sheaves
    \[\operatorname{Spl}_{\mathcal S/\Omega}(0,\mathcal H,1-g)\to \Omega\]
    with Mukai vector $(0,H_t,1-g)$ over $t\in\Omega$.
    With the choice of a generic stability condition on the central fiber, $M_{S_0}(0,H_0,1-g)$ is a compact subspace of $\operatorname{Spl}_{\mathcal S/\Omega}(0,\mathcal H,1-g)$ such that any two points of it are separated by open subsets.
    Therefore it has a separated neighbourhood which is proper over an open in $\Omega$ containing $0$.
\end{remark}

Write $|\mathcal H|\to\Omega$ for the relative linear system.
\begin{proposition}\label{prop:fibrationinfamily}
    There are maps
    \[\begin{tikzcd}
        M_{\mathcal S/\Omega}(0,\mathcal H,1-g)\rar[shift left]{\varpi}& {|\mathcal H|}\lar[shift left]{\xi}
    \end{tikzcd}\]
    so that $\varpi_t$ is the support morphism $M_{S_t}(0,H_t,1-g)\to |H_t|$, and $\xi_t$ is the zero-section $|H_t|\to M_{S_t}(0,H_t,1-g),C\mapsto[\mathcal O_C]$.
\end{proposition}
\begin{proof}
    It's clear that both maps exist for the corresponding moduli spaces of simple sheaves.
    Taking the stable moduli space does not prevent $\varpi$ from being a morphism.
    As for the section $\xi$, we need only to show that $\mathcal O_C$ is stable for our choice of polarization for all $C\in |H_0|$.

    Suppose $C=s+F_1+\cdots+F_g$.
    Write $H=s+kf$ for the generic polarization where $k$ is large.
    Let us compute
    \[P_H(\mathcal O_C)(n)=\frac{(nH)^2-(nH-C)^2}{2}=nH\cdot C-\frac{C^2}{2}=n(k+g-2)-(g-1).\]
    Suppose we have a destabilizing subsheaf of $\mathcal O_C$.
    Then it must be the ideal of a closed subscheme $D$ of $C$.
    Without loss of generality, we may assume that $D$ is pure of dimension $1$.

    If $D=s+F_1+\cdots+F_r$ for some $r$, then
    \begin{align*}
        P_H(\mathcal I_D)(n)&=P_H(\mathcal O_C)(n)-P_H(\mathcal O_D)(n)\\
        &=n(k+g-2)-(g-1)-n(k+r-2)+(r-1)\\
        &=n(g-r)-(g-r)
    \end{align*}
    which means that $\mathcal I_D$ is not destabilizing as $(g-1)/(k+g-2)<1$.
    Similarly, if $D=F_1+\cdots+F_r$, then
    \begin{align*}
        P_H(\mathcal I_D)(n)&=P_H(\mathcal O_C)(n)-P_H(\mathcal O_D)(n)\\
        &=n(k+g-2)-(g-1)-nH\cdot D=n(k+g-2-r)-(g-1)
    \end{align*}
    which again indicates that it is not destabilizing.
\end{proof}
\begin{remark}
    Under Yoshioka's deformation, $\xi_0$ corresponds to $s^{[g]}$.
\end{remark}
From now on let us fix such a deformation.
Write $\mathcal M=M_{\mathcal S/\Omega}(0,\mathcal H,1-g)$ and $M_t$ its fiber over $t\in\Omega$.
\subsection{Extended Mukai vectors}
To simplify the study of the Hodge theory of such a deformation, we use Beckmann's extended Mukai vectors to put everything inside the lattice $\tilde{H}(X;\mathbb Q)$.
Here, $X$ is a hyper-K\"ahler manifold of dimension $2d$.

For a class $v\in H^\ast(X;\mathbb Q)$, we write $\overline{v}$ for its projection (under the intersection pairing) to $\operatorname{SH}(X;\mathbb Q)$.
Similarly, we write
\[T:\operatorname{Sym}^d\tilde{H}(X;\mathbb Q)\to\operatorname{SH}(X;\mathbb Q)\]
for the projection to the image of $\Psi$, where the pairing on $\operatorname{Sym}^k\tilde{H}(X;\mathbb Q)$ is given by
\[b_{[k]}(\gamma_1\cdots\gamma_k,\gamma_1'\cdots\gamma_k')=(-1)^k\sum_{\sigma\in\operatorname{Sym}_k}\prod_{j=1}^kb(\gamma_j,\gamma_{\sigma(j)}').\]
Note that this pairing is nondegenerate.
Let us note that the pairing, hence the projection $T$, respects Hodge structures.
\begin{remark}\label{rem:Psi_isometry}
    When normalized correctly (\cite[Proposition 3.8]{taelman}), the injection $\Psi$ is an isometry if one takes the pairing on $\operatorname{SH}(X;\mathbb Q)$ to be the Mukai pairing
    \[b_{\operatorname{SH}}(\gamma_1\cdots\gamma_k,\gamma_1'\cdots\gamma_{2d-k}')=(-1)^k\int\gamma_1\cdots\gamma_k\gamma_1'\cdots\gamma_{2d-k}'.\]
\end{remark}
\begin{definition}
    Let $X$ be a hyper-K\"ahler variety.
    We say $v\in H^\ast(X;\mathbb Q)$ admits an extended Mukai vector if there is a nonzero $u\in \tilde{H}(X;\mathbb Q)$ and a nonzero constant $c\in\mathbb Q$ such that
    \[\overline{v}=cT(u^d).\]
\end{definition}
We call such a $u$ an extended Mukai vector for $v$.

The requirement of $c,u$ being nonzero enables us to obtain a uniqueness statement of $u$, which is also asserted in \cite{beckmann}.
It is not essential for our purpose.
We include here a statement and a proof for exposition.
\begin{proposition}
    If $v$ admits an extended Mukai vector, then it is unique up to a rational constant.
\end{proposition}
\begin{proof}
    Suppose we have
    \[cT(u^d)=c'T(u'^d)\]
    for some nonzero $c,c',u,u'$.
    We aim to show that $u=ru'$ for some $r\in\mathbb Q$.
    It suffices to show that $u=ru'$ for some $r\in\mathbb C$, since $u,u'$ are both rational vectors.
    So we may assume $c=c'=1$ without loss of generality.
    This means that $u^d-u'^d\in\ker T=(\ker\Delta)^\perp$, where $\Delta$ is the Laplacian defined in \Cref{thm:laplacian}.

    Let us compute this kernel after extending scalars to $\mathbb C$.
    Take an orthonormal basis $z_1,\ldots,z_s$ for $\tilde{H}(X;\mathbb C)$ under $b_{\mathbb C}$.
    Note that $s>2$.
    We can then identify $\operatorname{Sym}^k\tilde{H}(X;\mathbb C)$ with the degree $k$ homogeneous part of $\mathbb C[z_1,\ldots,z_s]$.

    Under this identification, it is easy to see that for any homogeneous $F,G$ of degree $k$,
    \[b_{[k]}(F,G)=(-1)^kF(\partial_{z_1},\ldots,\partial_{z_s})(G)|_{z=0}.\]
    Similarly, the Laplacian $\Delta:\operatorname{Sym}^d\tilde{H}(X;\mathbb C)\to \operatorname{Sym}^{d-2}\tilde{H}(X;\mathbb C)$ can be written as
    \[\Delta=\frac{1}{2}(\partial_{z_1}^2+\cdots+\partial_{z_s}^2).\]
    Consequently,
    \[b_{[d]}((z_1^2+\cdots+z_s^2)F,G)=(-1)^d2F(\partial_{z_1},\ldots,\partial_{z_s})\Delta G|_{z=0}=2b_{[d-2]}(F,\Delta G)\]
    for all $F\in \operatorname{Sym}^{d-2}\tilde{H}(X;\mathbb C),G\in\operatorname{Sym}^d\tilde{H}(X;\mathbb C)$.

    Since $b_{[k]}$ is always nondegenerate, this means that
    \[\ker T=(\ker\Delta)^\perp=(z_1^2+\cdots+z_s^2)\operatorname{Sym}^{d-2}\tilde{H}(X;\mathbb C).
    \]
    Hence $(z_1^2+\cdots+z_s^2)\mid u^d-u'^d$.
    But $s>2$, so $z_1^2+\cdots+z_s^2$ is irreducible, and yet $u^d-u'^d$ factorizes into linear polynomials $u-\zeta_d^ju'$.
    So there must be some $j$ such that $u-\zeta_d^ju'=0$, proving the claim.
\end{proof}
The construction of extended Mukai vectors also respects Hodge structures.
\begin{lemma}\label{lem:hodgegiveshodge}
    Suppose $v$ is of Hodge type and admits an extended Mukai vector $u$.
    Then $u\in\tilde{H}^{1,1}(X;\mathbb Q)$.
\end{lemma}
\begin{proof}
    If $v$ is of Hodge type, so is $\overline{v}$.
    Let us show that the only way for $T(u^d)$ to be of Hodge type is if $u$ is.
    Indeed, $\Psi(\sigma^d)=\sigma^d$, so if $\overline{v}$ is of Hodge type then, by \Cref{rem:Psi_isometry},
    \[0=b_{\operatorname{SH}}(\overline{v},\sigma^d)=b_{[d]}(T(u^d),\sigma^d)=b_{[d]}(u^d,\sigma^d)=b(u,\sigma)^d\]
    where each equality is up to multiplication by a nonzero constant.

    Hence $b(u,\sigma)=0$.
    Similarly $b(u,\bar\sigma)=0$.
    So $u\in\tilde{H}^{1,1}(X;\mathbb Q)$.
\end{proof}
Mukai vectors of various (but not all!) objects in the derived category admit extended Mukai vectors.
When they do, we have
\begin{theorem}[\cite{beckmann}]
    Suppose $\Phi:D^b\operatorname{Coh}(X)\to D^b\operatorname{Coh}(X)$ is a derived equivalence taking an object $A$ to $B$.
    If $v(A)$ admits an extended Mukai vector $u$, then $\Phi^{\tilde{H}}(u)$ is an extended Mukai vector for $v(B)$.
\end{theorem}
We are interested in the extended Mukai vectors for the section and fibers of the fibration $\pi^{[g]}:S_0^{[g]}\to\mathbb P^g$.

Let the classes $s,f,\delta,\alpha,\beta\in \tilde{H}(S_0^{[g]},\mathbb Q)$ be as before.
\begin{lemma}
    Let $F$ be a smooth fiber of $\pi^{[g]}$ and $s^{[g]}$ be the section as before.
    Then $v(\mathcal O_F)$ admits an extended Mukai vector
    \[\tilde{v}(\mathcal O_F)=f,\]
    and $v(\mathcal O_{s^{[g]}})$ admits an extended Mukai vector
    \[\tilde{v}(\mathcal O_{s^{[g]}})=s-\frac{\delta}{2}+\frac{g+1}{2}\beta.\]
\end{lemma}
Note that there is a sign error in \cite[Example 4.16]{beckmann} which is fixed here.
\begin{proof}
    The former follows from the same argument as in \cite[Example 4.17]{beckmann}.
    As for the latter, \cite[Remark 7.3]{beckmann} gives
    \[s-\frac{\delta}{2}-\frac{g+1}{2}\beta\]
    for the extended Mukai vector of $\omega_{s^{[g]}}=\mathcal O_{s^{[g]}}\otimes\pi^{[g]\ast}\mathcal O(-(g+1))$.
    So an extended Mukai vector for $\mathcal O_{s^{[g]}}$ is
    \[B_{(g+1)f}\left(s-\frac{\delta}{2}-\frac{g+1}{2}\beta\right)=s-\frac{\delta}{2}+\frac{g+1}{2}\beta\]
    as desired.
\end{proof}

\begin{remark}
    One may also make the same computation by computing the rest of the action of $(\mathcal F^{[g]})^{\tilde{H}}$ and the identities $\mathcal F^{[g]}(k_{\rm pt})=\mathcal O_F,\mathcal F^{[g]}(\mathcal O_{s^{[g]}})=\mathcal O_{S_0^{[g]}}$.
\end{remark}
The upshot is we can now identify a subspace of $\tilde{H}$ that remains of type $(1,1)$ under the deformation $\mathcal M\to\Omega$ considered before.
\begin{corollary}
    Let $V$ be as in \Cref{prop:actasscalaronperp}.
    Then $V$ remains of type $(1,1)$ under the parallel transport induced by the deformation
    \[M_{\mathcal S/\Omega}(0,\mathcal H,1-g)\to \Omega.\]
\end{corollary}
\begin{proof}
    The classes $\alpha,\beta$ remain of type $(1,1)$ by definition.
    By \Cref{prop:fibrationinfamily}, $v(\mathcal O_{s^{[g]}})$ gets transported to $v(\mathcal O_{\xi_t})$ which is algebraic, hence Hodge.
    So the extended Mukai vector of the former gets transported to the extended Mukai vector of the latter, which by \Cref{lem:hodgegiveshodge} is of type $(1,1)$.
    Therefore the transport of
    \[s-\frac{\delta}{2}+\frac{g+1}{2}\beta\]
    is of type $(1,1)$.
    Similarly, since the Mukai vector of (the structure sheaf of) a general fiber gets transported to the Mukai vector of a general fiber, $f$ also remains of type $(1,1)$ by the same argument.

    These classes generate $V$.
\end{proof}
\begin{corollary}\label{cor:remainhodge}
    Under the parallel transport induced by the same deformation, the Chern classes of $\mathcal P^{[g]}$ remain of Hodge type.
\end{corollary}
\begin{proof}
    It suffices to show that the Mukai vector of $\mathcal P^{[g]}$ remains of Hodge type under the deformation.

    We make use of \Cref{prop:actasscalaronperp}.
    The transported action $f$ of the Mukai vector of $\mathcal P^{[g]}$ still preserves the LLV algebra as well as the Verbitsky component, for they are topological invariants.
    Hence, it descends to an action $f^{\tilde{H}}$ on the extended Mukai lattice.

    However, since the transport of $V$ remains of type $(1,1)$ on this lattice, and $\Phi_{\mathcal P^{[g]}}^{\tilde{H}}$ acts as a constant on $V^\perp$, this transported action $f^{\tilde{H}}$ must preserve the Hodge structure.
    In other words, we have
    \[[f^{\tilde{H}},g]=0\]
    where $g$ is the Hochschild grading, i.e.~$g|_{H^{p,q}}=(p-q)\operatorname{id}$.
    Note that $g$ is an element of the complexified LLV algebra by \cite[Proposition 2.8]{taelman}.

    By \Cref{cor:commutewllv}, we then have $[g,f]=0$, i.e.~$f$ preserves the Hochschild grading on cohomology.
    Since $f$ is induced by the transport of the Mukai vector of $\mathcal P^{[g]}$, the latter remains of Hodge type.
\end{proof}
\begin{remark}\label{rem:yoshiokauniv}
    For later convenience, let us describe Yoshioka's deformation via periods.
    Write
    \[\mathcal D=\{\sigma\in\mathbb P(H^2(S_0,\mathbb C)):\sigma^2=0,\sigma\bar\sigma>0\},\]
    \[\mathcal D^{[g]}=\{\sigma\in\mathbb P(H^2(S_0^{[g]};\mathbb C)):\sigma^2=0,\sigma\bar\sigma>0\}\]
    for the period domains.
    We have then an isomorphism
    \[\{\sigma\in\mathcal D:\sigma\cdot H_0=0\}\to\left\{\sigma\in\mathcal D^{[g]}:\sigma\cdot f=\sigma\cdot\left(s-\frac{\delta}{2}\right)=0\right\}\]
    given by the period of the construction
    \[(S,H)\mapsto M_S(0,H,1-g)\]
    if the latter carries the marking
    $(0,H,1-g)_S^\perp\cong (0,s+gf,1-g)_{S_0}^\perp\cong H^2(S_0^{[g]};\mathbb Z)$ with the second isomorphism given by \Cref{thm:yoshiokaiso}.
    Explicitly, it acts as $-1$ on $\langle s,f\rangle^\perp=\langle s,f,\delta\rangle^\perp$ and sends $-s+(g-2)f$ to $\delta-(g-1)f$ by \Cref{prop:elliptic_surface_calc}.

    Therefore, a very general deformation of the fibration $S_0^{[g]}\to\mathbb P^g$ that preserves the section has to be an instance of Yoshioka's deformation.
\end{remark}
\subsection{The semiregularity theorem}
In order to deform the BKR congruence $\mathcal P^{[g]}$ away from the central fiber, we make use of the semiregularity theorem due to Buchweitz and Flenner \cite{semireg}.
\begin{definition}
    Let $\mathcal G$ be a coherent sheaf on a smooth proper variety $M$.
    Its semiregularity map is
    \[\operatorname{Ext}^2(\mathcal G,\mathcal G)\to\bigoplus_pH^{2+p}(M,\Omega_M^p)=\operatorname{HH}_{-2}(M)\]
    where each $\operatorname{Ext}^2(\mathcal G,\mathcal G)\to H^{2+p}(M,\Omega_M^p)$ is given by the composition
    \[\begin{tikzcd}
        \operatorname{Ext}^2(\mathcal G,\mathcal G)\rar{\operatorname{At}_p(\mathcal G)/p!}&\operatorname{Ext}^{2+p}(\mathcal G,\mathcal G\otimes\Omega_M^p)\rar{\operatorname{tr}}& H^{2+p}(M,\Omega_M^p)
    \end{tikzcd},\]
    where $\operatorname{At}_p$ is the $p$-th Atiyah class.

    We say $\mathcal G$ is semiregular if its semiregularity map is injective.
\end{definition}
Informally, the semiregularity map takes the obstruction to deforming $\mathcal G$ to the obstruction to keeping $\operatorname{ch}\mathcal G$ as a Hodge class in the deformation.
One precise statement is the following:
\begin{theorem}[{\cite[Theorem 5.1]{semireg}}]
    Let $\Omega\ni 0$ be a small disk over which we have a smooth proper family $\mathcal M=(M_t)_{t\in\Omega}\to\Omega$.
    Suppose $\mathcal G_0\in\operatorname{Coh}(M_0)$ is semiregular, and the parallel transports of $\operatorname{ch}\mathcal G_0$ remain Hodge, then the local deformation of $\mathcal G_0$ in this family is unobstructed, i.e.~it admits a flat deformation $\mathcal G_t$ for $t$ near $0$.
\end{theorem}
In general, it is quite rare that one can find a semiregular sheaf.
Indeed, $\operatorname{HH}_{-2}(M)$ may well be smaller than $\operatorname{Ext}^2(\mathcal G,\mathcal G)$.
For Fourier-Mukai kernels, however, this $\operatorname{Ext}$ group can be better controlled.
In fact, we have:
\begin{theorem}[{\cite[Theorem 8.4]{rienks}}]
    Suppose $\mathcal G$ is the kernel of a fully faithful Fourier-Mukai transform
    \[\Phi_{\mathcal G}:D^b\operatorname{Coh}(X)\to D^b\operatorname{Coh}(Y),\]
    with $\omega_X\cong\mathcal O_X$, then $\mathcal G$ is semiregular.
\end{theorem}
Let us apply this to the kernel $\mathcal P^{[g]}$ deforming in the family $\mathcal M\times_\Omega\mathcal M\to\Omega$.
\begin{proposition}
    Possibly after shrinking $\Omega$, the kernel $\mathcal G_0=\mathcal P^{[g]}$ admits a locally unique flat deformation $\mathcal G_t$ on $M_t\times M_t$ such that each $\mathcal G_t$ is Cohen-Macaulay, and defines a Fourier-Mukai autoequivalence of $D^b\operatorname{Coh}(M_t)$.
\end{proposition}
\begin{proof}
    The preceding theorem shows that $\mathcal G_0$ is semiregular.
    Moreover, its Chern character remains Hodge in this deformation by \Cref{cor:remainhodge}.
    The theorem of Buchweitz and Flenner applies to give a flat deformation $\mathcal G_t$.

    Whether a sheaf is Cohen-Macaulay, and whether it induces a Fourier-Mukai equivalence are both open conditions, so $\mathcal G_t$ is an equivalence possibly after shrinking $\Omega$.

    The local uniqueness of this deformation follows from the rigidity of $\mathcal G_0$, as supplied by the next lemma.
\end{proof}

\begin{lemma}\label{lem:FMkernelsarerigid}
    Suppose $X,Y$ are smooth projective varieties with $\operatorname{HH}^1(X)=0$.
    Then any kernel $\mathcal K$ of a derived equivalence $D^b\operatorname{Coh}(X)\to D^b\operatorname{Coh}(Y)$ is rigid.
\end{lemma}
\begin{proof}
    Composing with the inverse kernel gives
    \[\operatorname{Ext}^1(\mathcal K,\mathcal K)=\operatorname{Ext}^1(\mathcal O_\Delta,\mathcal O_\Delta)=\operatorname{HH}^1(X)=0,\]
    as in the proof of \cite[Theorem 8.4]{rienks}.
\end{proof}

Let us note some properties of this deformation.
\begin{proposition}\label{lem:propofdef}
    Possibly after shrinking $\Omega$, $\mathcal G_t$ is invariant under swapping the factors, and satisfies
    \[\mathcal G_t\otimes\operatorname{pr}_1^\ast\varpi_t^\ast\mathcal O(1)\cong \mathcal G_t\otimes\operatorname{pr}_2^\ast\varpi_t^\ast\mathcal O(1).\]
\end{proposition}

\begin{proof}
    The sheaf $\mathcal G_0^\top$ (where ${}^\top$ denotes swapping the factors) has a flat deformation, namely $\mathcal G_t^\top$.
    But $\mathcal G_0\cong\mathcal G_0^\top$, so uniqueness gives $\mathcal G_t\cong\mathcal G_t^\top$.

    Similarly, $\mathcal G_0\otimes \operatorname{pr}_1^\ast(\pi^{[g]})^\ast\mathcal O(1)\otimes \operatorname{pr}_2^\ast(\pi^{[g]})^\ast\mathcal O(-1)$ deforms to
    \[\mathcal G_t\otimes \operatorname{pr}_1^\ast\varpi_t^\ast\mathcal O(1)\otimes \operatorname{pr}_2^\ast\varpi_t^\ast\mathcal O(-1).\]
    But we have $\mathcal G_0\otimes \operatorname{pr}_1^\ast(\pi^{[g]})^\ast\mathcal O(1)\cong\mathcal G_0\otimes \operatorname{pr}_2^\ast(\pi^{[g]})^\ast\mathcal O(1)$ since $\mathcal G_0$ is supported on $S_0^{[g]}\times_{\mathbb P^g}S_0^{[g]}$ (see \cite[Lemma 1.5]{bai25}).
    Uniqueness again shows
    \[\mathcal G_t\otimes\operatorname{pr}_1^\ast\varpi_t^\ast\mathcal O(1)\cong \mathcal G_t\otimes\operatorname{pr}_2^\ast\varpi_t^\ast\mathcal O(1).\qedhere\]
\end{proof}

\subsection{Support of the deformation}
Since we are deforming from a relative transform, we certainly expect the deformation to be a relative transform as well.
The idea is to use the fact that $\Phi_{\mathcal G_t}$ commutes with $-\otimes\varpi_t^\ast\mathcal O(1)$, which follows from the isomorphism
\[\mathcal G_t\otimes\operatorname{pr}_1^\ast\varpi_t^\ast\mathcal O(1)\cong \mathcal G_t\otimes\operatorname{pr}_2^\ast\varpi_t^\ast\mathcal O(1)\]
in \Cref{lem:propofdef}.
This is necessary for something to be a relative kernel \cite[Lemma 1.5]{bai25}.
The idea is to show that it is more or less sufficient by, essentially, generating the rest of $D^b\operatorname{Coh}(|H_t|)$ from $\mathcal O(1)$.

This doesn't quite work as stated, since it does not detect an automorphism of the base $|H_t|$.
But it can be remedied in our situation:
\begin{lemma}\label{aicorrectedthis}
    Suppose there is an automorphism $f$ of $|H_t|$ such that the scheme-theoretic support of $\mathcal G_t$ is contained in the $f$-twisted fiber product, that is the preimage of the graph of the automorphism under
    \[M_t\times M_t\to |H_t|\times |H_t|,\]
    then the scheme-theoretic support of $\mathcal G_t$ is precisely this twisted fiber product, and $f^2=\operatorname{id}$.
\end{lemma}
\begin{proof}
    First of all, since the $f$-twisted fiber product is integral and $\mathcal G_t$ is a flat deformation of $\mathcal G_0$, the scheme-theoretic support is precisely the $f$-twisted fiber product.

    Now we also know that $\mathcal G_t^\top=\mathcal G_t$, hence the graph of this automorphism is stable under swapping the factors.
    But this can only happen if $f=f^{-1}$.
\end{proof}
\begin{lemma}
    The identity is an isolated point of the locus of involutions in $\operatorname{PGL}_{g+1}$.
\end{lemma}
\begin{proof}
    This is simply because the conjugacy classes are already closed, and any involution is diagonalizable.
\end{proof}
\begin{proposition}
    The scheme-theoretic support of $\mathcal G_t$ is $M_t\times_{|H_t|}M_t$.
\end{proposition}
\begin{proof}
    We make use of the isomorphism
    \[\tau:\mathcal G_t\otimes\operatorname{pr}_1^\ast\varpi_t^\ast\mathcal O(1)\cong \mathcal G_t\otimes\operatorname{pr}_2^\ast\varpi_t^\ast\mathcal O(1).\]

    Composing with the inverse kernel gives isomorphisms
    \begin{align*}
        \operatorname{Hom}(\mathcal G_t,\mathcal G_t\otimes\operatorname{pr}_i^\ast\varpi_t^\ast\mathcal O(1))&\cong\operatorname{Hom}(\mathcal O_\Delta,\mathcal O_\Delta\otimes\operatorname{pr}_i^\ast\varpi_t^\ast\mathcal O(1))\\
        &\cong H^0(M_t,\varpi_t^\ast\mathcal O(1))=H^0(|H_t|,\mathcal O(1))
    \end{align*}
    since $\varpi_{t\ast}\mathcal O_{M_t}=\mathcal O_{|H_t|}$.
    Written explicitly, a section $\mathcal O_{M_t}\to\varpi_t^\ast\mathcal O(1)$ simply corresponds to the homomorphism
    \[\mathcal G_t=\mathcal G_t\otimes \operatorname{pr}_i^\ast\mathcal O_{M_t}\to \mathcal G_t\otimes\operatorname{pr}_i^\ast\varpi_t^\ast\mathcal O(1).\]

    The isomorphism $\tau$ then gives rise to an automorphism $\tilde{f}:H^0(M_t,\varpi_t^\ast\mathcal O(1))\to H^0(M_t,\varpi_t^\ast\mathcal O(1))$ which projectivizes into an automorphism $f=\mathbb P\tilde{f}:|H_t|\to |H_t|$.
    Let us show that the scheme-theoretic support of $\mathcal G_t$ is contained in the $f$-twisted fiber product.

    Let $s_0,s_1,\ldots,s_g$ be a basis for $H^0(M_t,\varpi_t^\ast\mathcal O(1))$ and let $t_j=\tilde{f}(s_j)$.
    Note that the $f$-twisted fiber product is cut out by the collection
    \[u_{ij}=s_i\boxtimes t_j-s_j\boxtimes t_i\in H^0(M_t\times M_t,\varpi_t^\ast\mathcal O(1)\boxtimes \varpi_t^\ast\mathcal O(1)).\]
    Let us show that they annihilate $\mathcal G_t$.
    This is clear: on opens in which both $\operatorname{pr}_i^\ast\varpi_t^\ast\mathcal O(1)$ trivialize, any section $v$ of $\mathcal G_t$ has $\tau(v\operatorname{pr}_1^\ast s_j)=v\operatorname{pr}_2^\ast t_j$, so if we write $\tau^+=\tau\boxtimes\operatorname{pr}_2^\ast\varpi_t^\ast\mathcal O(1)$ then
    \[\tau^+(vu_{ij})=\tau^+(v(s_i\boxtimes t_j-s_j\boxtimes t_i))=vt_i\boxtimes t_j-vt_j\boxtimes t_i=0,\]
    which means that $vu_{ij}=0$.

    There is therefore a continuous family $f_t\in\operatorname{PGL}_{g+1}$ such that the scheme-theoretic support of $\mathcal G_t$ is the $f_t$-twisted fiber product.
    Since $f_0=\operatorname{id}$, the preceding lemmas imply $f_t=\operatorname{id}$ for all $t$.
\end{proof}
\begin{remark}
    Alternatively, one may show this by comparing the deformation theory of the subvariety $M_0\times_{|H_0|}M_0$ in this family with $\mathfrak{pgl}_{g+1}$, at a cost of possibly shrinking $\Omega$ further.
\end{remark}
\begin{corollary}
    The sheaf $\mathcal G_t$ is a maximal Cohen-Macaulay sheaf on $M_t\times_{|H_t|}M_t$.
\end{corollary}
In particular, it is a line bundle over the smooth locus of $M_t\times_{|H_t|}M_t$.

Recall that the section $s^{[g]}:\mathbb P^g\to S_0^{[g]}$ deforms to a section $\xi_t:|H_t|\to M_t$.
The original kernel $\mathcal P^{[g]}$ is normalized at this section.
Naturally, one expects a similar normalization to hold for the deformation.
\begin{corollary}\label{cor:normalized_def}
    We have $\mathcal G_t|_{\xi_t\times_{|H_t|}M_t}\cong\mathcal O_{M_t}$ and $\mathcal G_t|_{M_t\times_{|H_t|}\xi_t}\cong\mathcal O_{M_t}$.
\end{corollary}
\begin{proof}
    This essentially follows from the simply-connectedness of $M_t$.
    Indeed, the sheaf $\mathcal G_t|_{\xi_t\times_{|H_t|}M_t}$ is a line bundle by the preceding corollary, and it is also a deformation of $\mathcal O_{M_0}$.
    That means $c_1(\mathcal G_t|_{\xi_t\times_{|H_t|}M_t})=0$, hence this line bundle is trivial as $H^1(M_t,\mathcal O_{M_t})=0$.

    The other side follows from $\mathcal G_t^\top=\mathcal G_t$.
\end{proof}

\section{Comparison with Arinkin's sheaf}\label{s:comparison}
\subsection{Generic fiber of the deformation}
Let $A$ be an abelian variety, and let $\mathcal L\in\operatorname{Pic}(A\times \check{A})$ be the Poincar\'e line bundle.
Suppose we are given an isomorphism $\theta:A\to\check{A}$ (e.g.~a principal polarization), then we get a line bundle $L=(\operatorname{id}\times\theta)^\ast\mathcal L$ which induces an autoequivalence of $D^b\operatorname{Coh}(A)$.

The ambiguity introduced by $\theta$ can be recovered from the cohomological action of this Fourier-Mukai autoequivalence.
The associated cohomological action $\Phi_L^H$ is
\[\alpha\mapsto\operatorname{pr}_{2\ast}(\exp(c_1(L))\smile \operatorname{pr}_1^\ast\alpha).\]

Note that $\exp(c_1(L))$ is integral, so this actually gives a map $H^\ast(A;\mathbb Z)\to H^\ast(A;\mathbb Z)$.
As $L$ is the pullback of the Poincar\'e bundle, it is also degree-reversing.
In particular, it restricts to a map $\theta_L:H^{2g-1}(A;\mathbb Z)\to H^1(A;\mathbb Z)$.
\begin{lemma}\label{lem:normalizedlb_c1recovery}%... look over again
    Suppose $A$ is an abelian variety of dimension $g$ and $L\in\operatorname{Pic}(A\times A)$ is the kernel of a Fourier-Mukai autoequivalence of $A$ with
    \[L|_{\{0_A\}\times A}=\mathcal O_A,L|_{A\times\{0_A\}}=\mathcal O_A.\]

    Then $L$ can be recovered from $c_1(L)$, or equivalently, the action
    \[\theta_L:H^{2g-1}(A;\mathbb Z)\to H^1(A;\mathbb Z).\]
\end{lemma}
In fact, $\theta_L$ is precisely the classifying isomorphism $\theta:A\to\check{A}$ so that $L\cong(\operatorname{id}\times\theta)^\ast\mathcal L$.
\begin{proof}
    Note that the class $c_1(L)$ and the map $\theta_L$ contain the same information by the K\"unneth formula.
    The set of all line bundles on $A\times A$ with the same $c_1$ are exactly $L\otimes(M\boxtimes N)$ for $M,N\in\operatorname{Pic}^0(A)$.

    But by the normalizations we have
    \[L\otimes(M\boxtimes N)|_{\{0_A\}\times A}=N,L\otimes(M\boxtimes N)|_{A\times \{0_A\}}=M,\]
    so they cannot be $L$ unless $M=N=\mathcal O$.
\end{proof}
The point of all these is that the polarization is a piece of discrete data.

Now let $C$ be a curve.
We write $\operatorname{Jac}(C)$ for the connected component of $\operatorname{Pic}(C)$ containing the identity.
Suppose in addition that $C$ is of compact type.
Then the canonical principal polarization gives rise to a line bundle $L$ over $\operatorname{Jac}(C)\times\operatorname{Jac}(C)$.
If one chooses a theta divisor $\Theta$ on $\operatorname{Jac}(C)$, then we can describe
\[L\cong m^\ast\Theta\otimes\operatorname{pr}_1^\ast\Theta^\vee\otimes\operatorname{pr}_2^\ast\Theta^\vee\]
where $m:\operatorname{Jac}(C)\times\operatorname{Jac}(C)\to\operatorname{Jac}(C)$ is the multiplication.
We call $L$ the canonical Poincar\'e bundle for this Jacobian.
\begin{proposition}
    Let $C$ be a general member of $|s+gf|$ on $S_0$.
    In particular, it has compact type.
    Then
    \[\mathcal P^{[g]}|_{[C]}\cong L\]
    where $L$ is the canonical Poincar\'e bundle for $\operatorname{Jac}(C)$.
\end{proposition}
Here, the Jacobian of any $C\in |H_t|$ of compact type is identified with the fiber of $\varpi_t$ at $[C]$ by sending $M\in\operatorname{Jac}(C)$ to $i_{C\ast}M\in M_{S_t}(0,H_t,1-g)$.
This lands in the stable locus because $\mathcal O_C$ does.
\begin{proof}
    Write the elliptic fibers in $C$ as $E_1,\ldots,E_g$.
    The generality assumption says $E_j$ are distinct and smooth.
    We can then identify $\operatorname{Jac}(C)\cong E_1\times\cdots\times E_g$, and
    \[\Phi_L(k_{(p_1,\ldots,p_g)})=\mathcal O_{E_1}(0_{E_1}-p_1)\boxtimes\cdots\boxtimes\mathcal O_{E_g}(0_{E_g}-p_g).\]

    On the other hand, let us write
    \[\begin{tikzcd}
        S_0^{[g]}&{[(\tilde{S}_0)_g/\operatorname{Sym}_g]}\rar\lar&{[S_0^g/\operatorname{Sym}_g]}
    \end{tikzcd}.\]
    Note that since $E_j$ are distinct fibers, any fiber of $\pi^{[g]}$ over $\pi(E_1)+\cdots+\pi(E_g)$ is contained in the locus in $S_0^{[g]}$ consisting of $g$ distinct points.
    View $p_1,\ldots,p_g$ as points on $S_0$, then
    \begin{align*}
        \mathcal F^{\boxtimes g}\operatorname{BKR}(k_{p_1+\cdots+p_g})&=\mathcal F^{\boxtimes g}(\operatorname{ind}_1^{\operatorname{Sym}_g}k_{(p_1,\ldots,p_g)})\\
        &=\operatorname{ind}_1^{\operatorname{Sym}_g}\left(i_{E_1\ast}\mathcal O_{E_1}(0_{E_1}-p_1)\boxtimes\cdots\boxtimes i_{E_g\ast}\mathcal O_{E_g}(0_{E_g}-p_g)\right),
    \end{align*}
    which $\operatorname{BKR}^\top$ then sends to
    \[i_{(\pi^{[g]})^{-1}(\pi(E_1)+\cdots+\pi(E_g))\ast}(\mathcal O_{E_1}(0_{E_1}-p_1)\boxtimes\cdots\boxtimes\mathcal O_{E_g}(0_{E_g}-p_g))\qedhere\]
\end{proof}
\begin{corollary}
    Suppose $C\in |H_t|$ is smooth, then $\mathcal G_t|_{[C]}$ is the canonical Poincar\'e bundle for $\operatorname{Jac}(C)$.
\end{corollary}
\begin{proof}
    Note that both sides are deformations of the situation in the preceding proposition.
    In particular, they have the same $c_1$.
    But we also know that both of them are normalized at the zero sections by \Cref{cor:normalized_def}.
    The result follows from \Cref{lem:normalizedlb_c1recovery}.
\end{proof}
\subsection{Arinkin's sheaf}
Let us now compare $\mathcal G_t$ to the construction of Arinkin \cite{arinkin}.
For an integral curve $C$ with planar singularities, we write $\overline{\operatorname{Jac}}(C)$ for its compactified Jacobian, parameterizing rank $1$ torsion-free sheaves on $C$.

Arinkin considered the Poincar\'e bundle $P$ over
\[U=(\operatorname{Jac}(C)\times\overline{\operatorname{Jac}}(C))\cup(\overline{\operatorname{Jac}}(C)\times\operatorname{Jac}(C))\subset \overline{\operatorname{Jac}}(C)\times\overline{\operatorname{Jac}}(C)\]
whose fiber at $(L,L')$ is
\[\det\mathbf R\Gamma(L\otimes L')\otimes \det\mathbf R\Gamma(\mathcal O_C)\otimes (\det\mathbf R\Gamma(L))^{-1}\otimes (\det\mathbf R\Gamma(L'))^{-1}.\]
Note that when $C$ is a smooth curve, then $P$ is precisely the dual of the kernel of the autoequivalence induced by the usual principal polarization.

Arinkin constructed a maximal Cohen-Macaulay extension $\bar{P}$ of $P$ over $\overline{\operatorname{Jac}}(C)\times\overline{\operatorname{Jac}}(C)$, and he showed that it induces an autoequivalence of $D^b\operatorname{QCoh}(\overline{\operatorname{Jac}}(C))$.
Note that such a maximal Cohen-Macaulay extension is necessarily unique, as the complement of $U$ has codimension $2$.

Now let $t\neq 0$.
The fiber of $\varpi_t:M_t\to |H_t|$ over any $C\in |H_t|$ (necessarily integral by the choice we have made) can be identified with $\overline{\operatorname{Jac}}(C)$ by sending $M\in \overline{\operatorname{Jac}}(C)$ to $i_{C\ast}M\in M_{S_t}(0,H_t,1-g)$.
Arinkin's construction is functorial, therefore gives rise to a maximal Cohen-Macaulay sheaf
\[\mathcal G_{\rm Arinkin}\in\operatorname{Coh}(M_t\times_{|H_t|}M_t)\]
inducing an autoequivalence of $D^b\operatorname{Coh}(M_t)$.

\begin{proposition}
    Suppose $\mathcal G$ is a maximal Cohen-Macaulay sheaf on $M_t\times_{|H_t|}M_t$ such that
    \begin{enumerate}
        \item $\mathcal G|_{\xi_t\times_{|H_t|}M_t}\cong\mathcal O_{M_t},\mathcal G|_{M_t\times_{|H_t|}\xi_t}\cong\mathcal O_{M_t}$, and
        \item the fiber of $\mathcal G$ over a smooth curve $C\in |H_t|$ is the canonical Poincar\'e bundle for $\operatorname{Jac}(C)$.
    \end{enumerate}
    Then $\mathcal G\cong \mathcal G_{\rm Arinkin}^\vee$.
\end{proposition}
Recall that the derived dual of a maximal Cohen-Macaulay sheaf is a maximal Cohen-Macaulay sheaf.
\begin{proof}
    Let us write $M^\circ_t\subset M_t$ for the locus that is the open $\operatorname{Jac}(C)\subset\overline{\operatorname{Jac}}(C)$ on each fiber.
    Write $U=(M_t^\circ\times_{|H_t|}M_t)\cup (M_t\times_{|H_t|}M_t^\circ)$.
    Since its complement has codimension at least $2$, it suffices to show an isomorphism of line bundles $\mathcal G|_U\cong\mathcal G_{\rm Arinkin}^\vee|_U$.

    Write $|H_t|^\circ\subset |H_t|$ for the locus of smooth curves.
    The second condition, combined with the seesaw theorem, shows that
    \[\mathcal G_{\rm Arinkin}^\vee|_{|H_t|^\circ}\cong\mathcal G|_{|H_t|^\circ}\otimes L\]
    for some line bundle $L$ pulled back from $|H_t|^\circ$.
    The normalizations imply that the restrictions of both $\mathcal G$ and $\mathcal G_{\rm Arinkin}^\vee$ to $\xi_t\times_{|H_t|}\xi_t$ are the trivial line bundle.
    Therefore $L$ is trivial.

    Take a section of $\mathcal G_{\rm Arinkin}|_U\otimes\mathcal G|_U$ trivializing it over $|H_t|^\circ$.
    Since $\varpi_t$ has integral fibers, so does the fibration $M_t\times_{|H_t|}M_t\to |H_t|$, hence $U\to |H_t|$.
    So the divisor of this section is the pullback of a divisor on $|H_t|$ supported in $|H_t|\setminus |H_t|^\circ$ via $U\to |H_t|$.
    In other words, $\mathcal G_{\rm Arinkin}|_U\otimes\mathcal G|_U$ is a line bundle pulled back from the base.
    The trivialization on $\xi_t\times_{|H_t|}\xi_t$ again shows that this line bundle is trivial.
\end{proof}
\begin{corollary}
    The deformation $\mathcal G_t$ is precisely the dual of $\mathcal G_{\rm Arinkin}$.
\end{corollary}
This deformation is thus another way to construct $\mathcal G_{\rm Arinkin}$.

\subsection{An application: the \texorpdfstring{$P=C$}{P=C} property}
As a first application of this degeneration, let us observe that the parallel transport of $\operatorname{ch}(\mathcal P^{[g]})$ would be precisely $\operatorname{ch}(\mathcal G_t)=\operatorname{ch}(\mathcal G_{\rm Arinkin}^\vee)$.
Consequently, suitably topological assertions about the Chern character of $\mathcal G_0=\mathcal P^{[g]}$ are equivalent to those same assertions about the Chern character of $\mathcal G_t$.

Let us give a concrete example of an application of this.
Recall that for smooth projective varieties $M,B$ and a fibration $h:M\to B$, the perverse filtration is a filtration on $H^\ast(M;\mathbb Q)$ given by
\[P_kH^l(M;\mathbb Q)=\operatorname{Im}(H^l(B,{}^{\mathfrak p}\tau_{\le k+\dim B}\mathbf Rh_\ast\underline{\mathbb Q})\to H^l(B,\mathbf Rh_\ast\underline{\mathbb Q})=H^l(M;\mathbb Q)),\]
where ${}^{\mathfrak p}\tau$ are perverse truncation functors.

This is an entirely topological notion.
It is, therefore, extremely interesting that the perverse filtration can sometimes be recovered algebraically.
This is the content of the $P=C$ property.
For us, it will mean the following phenomenon:

Let $\mathcal G\in D^b\operatorname{Coh}(M\times_B M)$ and write $i:M\times_BM\to M\times M$.
We write
\[\operatorname{ch}(i_\ast \mathcal G)\operatorname{td}(B)=\mathfrak F(\mathcal G)=\sum_{j\ge 0}\mathfrak F_j(\mathcal G),\mathfrak F_j(\mathcal G)\in\operatorname{Corr}_B^{j-g}(M,M),\]% also H^{2(j+\dim B)}(M\times M;\mathbb Q) if in cohomology
where $g=\dim M-\dim B$.

\begin{definition}
    We say that $M\to B$ satisfies the $P=C$ property realized by $\mathcal G$ if, when realized on cohomology, the correspondences $\mathfrak F_j(\mathcal G)$ satisfy
    \[\sum_{j=0}^k\operatorname{Im}\mathfrak F_j(\mathcal G)=P_kH^\ast(M;\mathbb Q).\]
\end{definition}
Roughly speaking, this says that the perverse filtration, which is a topological construction, can be reconstructed using algebraic operators.

An easy calculation shows that $\mathfrak F_j(\mathcal G^\vee)=(-1)^j\mathfrak F_j(\mathcal G)$ where the dual is taken in $M\times_BM$.
So $P=C$ is realized by $\mathcal G$ if and only if it is realized by $\mathcal G^\vee$.

The following proposition is then a direct consequence of the degeneration.
\begin{proposition}\label{prop:cohP=C}
    Let $(S,H)$ be a polarized K3 surface with $|H|$ integral and of genus $g$.
    Let $S_0\to\mathbb P^1$ be an elliptic K3 surface with integral fibers and a section.
    Then $M_S(0,H,1-g)\to |H|$ satisfies $P=C$ realized by Arinkin's sheaf if and only if $S_0^{[g]}\to\mathbb P^g$ satisfies $P=C$ realized by $\mathcal P^{[g]}$.
\end{proposition}
\begin{proof}
    Note that the perverse filtration is stable under parallel transport along a family of fibrations by \cite[Theorem 3.2.1]{perverselc}.
    Applying this to the Yoshioka deformation shows that $P=C$ for $\mathcal P^{[g]}$ is equivalent to $P=C$ for Arinkin's sheaf in a nearby Beauville-Mukai system.

    The space of $(S,H)$ with $|H|$ integral is connected, so the same degeneration argument applies to get to the prescribed Beauville-Mukai system.
\end{proof}
In \cite{msy}, the authors developed a general framework to study what they call \emph{dualizable abelian fibrations}.
In particular, they showed that the Beauville-Mukai system $M_S(0,H,1-g)\to |H|$ satisfies $P=C$ realized by Arinkin's sheaf (\cite[Theorem 0.2, Theorem 2.6]{msy}).
We therefore have:
\begin{corollary}
    The Hilbert scheme fibration $S_0^{[g]}\to\mathbb P^g$ satisfies $P=C$ realized by $\mathcal P^{[g]}$.
\end{corollary}
In upcoming work of Ghosh \cite{soumik}, this corollary is proven independently using Nakajima operators.
Using his work, this degeneration then gives a different proof of $P=C$ for general Beauville-Mukai systems.

One can also obtain a categorical version of \Cref{prop:cohP=C} similar to \cite[Theorem 2.6]{msy}.
Suppose $\mathcal G$ induces a Fourier-Mukai autoequivalence.
Denote by
\[\mathfrak F^{-1}(\mathcal G)=\sum_{j\ge 0}\mathfrak F_j^{-1}(\mathcal G),\mathfrak F_j^{-1}(\mathcal G)\in\operatorname{Corr}_B^{j-g}(M,M)\]
the inverse correspondence.
\begin{definition}[{\cite[Definition 2.5]{msy}}]
    We say $\mathcal G$ satisfies Fourier vanishing if
    \[\mathfrak F^{-1}_a(\mathcal G)\circ\mathfrak F_b(\mathcal G)=0,\]
    for any $a+b<2g$.
\end{definition}
Write
\[\mathfrak p_k(\mathcal G)=\sum_{j\le k}\mathfrak F_j(\mathcal G)\circ\mathfrak F^{-1}_{2g-j}(\mathcal G),\]
following \cite[Theorem 2.6]{msy}.
Note that if $\mathcal G$ satisfies Fourier vanishing, then $\mathfrak p_k$ is idempotent, i.e.~a projector.

By \cite[\S3.5]{msy}, Arinkin's sheaf satisfies Fourier vanishing.
Our deformation immediately implies:
\begin{corollary}
    The sheaf $\mathcal G=\mathcal P^{[g]}$ satisfies Fourier vanishing.
\end{corollary}
The cycles $\mathfrak p_k$ may be understood as elements of
\[\operatorname{End}(\mathbf Rh_\ast\underline{\mathbb Q})\cong H^{\rm BM}_{2\dim M}(M\times_BM;\mathbb Q)\]
via the Corti-Hanamura realization functor \cite{corti_hanamura}.
\begin{definition}[Categorical $P=C$]
    We say $\mathcal G$ realizes categorical $P=C$ for $h:M\to B$ if, under the Corti-Hanamura realization, the map
    \[\operatorname{im}\mathfrak p_k(\mathcal G)\hookrightarrow \mathbf Rh_\ast\underline{\mathbb Q}\]
    is precisely the perverse truncation
    \[{}^{\mathfrak p}\tau_{\le k+\dim B}\mathbf Rh_\ast\underline{\mathbb Q}\hookrightarrow \mathbf Rh_\ast\underline{\mathbb Q}.\]
\end{definition}
General Beauville-Mukai systems satisfy this categorical $P=C$ property by \cite[Theorem 0.2, Theorem 2.6]{msy}.

Let us now place ourselves in the situation considered in \Cref{prop:cohP=C}.
We then have:
\begin{proposition}\label{prop:P=Ccat}
    A general Beauville-Mukai system satisfies categorical $P=C$ if and only if the Hilbert scheme fibration does.
\end{proposition}
\begin{proof}
    Since $\mathfrak p_k$ is a projector, the $P=C$ statement can alternatively be written as saying
    \[{}^{\mathfrak p}\mathcal H^i(\mathfrak p_k(\mathcal G))=\operatorname{id}\]
    for $i\le k+\dim B$, and zero otherwise.

    Take the Yoshioka deformation $\mathcal M\to |\mathcal H|\to\Omega$ like before.
    Then $\phi_t[-1]\mathbf R\varpi_\ast\underline{\mathbb Q}=0$, where $\phi_t[-1]$ is the vanishing cycles over $t\in\Omega$.
    Write $P_i={}^{\mathfrak p}\mathcal H^i(\mathbf R\varpi_\ast\underline{\mathbb Q})$.
    Then the map
    \[\operatorname{End}(P_i)\to\operatorname{End}(i_t^\ast[-1] P_i)\]
    is injective.
    Indeed, suppose $f:P_i\to P_i$ is any map whose restriction is zero, and we take any simple factor $Q$ of its image.
    Then $\phi_s[-1]Q=0$ for all $s$.
    In particular, the support of $Q$ must surject onto $\Omega$.
    But then the nearby cycles functor $\psi_s[-1]Q$ can never vanish.
    Hence $i_t^\ast[-1]Q=\psi_t[-1]Q\neq 0$.

    One can therefore lift the $P=C$ statement at a specific fiber to every fiber.
\end{proof}

\section{A factorization of Markman's bundle}\label{s:hyperhol}

\subsection{Theorem of the square and twists}
For $t\neq 0$, the fibration $\varpi_t:M_t\to |H_t|$ is a compactified Jacobian fibration for the universal curve over $|H_t|$.
So if we take $M_t^\circ\subset M_t$ the relative Jacobian of the curves, then we have a morphism
\[\mu_t:M_t^\circ\times_{|H_t|}M_t\to M_t,(L,L')\mapsto L\otimes L'.\]

On $M_0$, one can similarly consider the locus $M_0^\circ\subset M_0$ as described in \cite{bai25}.
When interpreted as an open inside the moduli $M_0=M_{S_0}(0,s+gf,1-g)$, $M_0^\circ$ parameterizes (pushforwards) of line bundles that have degree $0$ on the reduction of each irreducible component of a curve.
As explained in \cite{bai25}, we also have a morphism
\[\mu_0:M_0^\circ\times_{|H_0|}M_0\to M_0,(L,L')\mapsto L\otimes L'.\]

Denote by $\mathcal M^\circ\subset\mathcal M$ the locus whose fiber over $t\in\Omega$ is $M_t^\circ$.
Then these $\mu_t$ can be collected into a morphism
\[\mu:\mathcal M^\circ\times_{|\mathcal H|}\mathcal M\to\mathcal M,\]
and it is a group action.
More precisely,
\begin{lemma}
    The morphism $\mu$ restricts to a commutative group law making $\mathcal M^\circ$ an $|\mathcal H|$-group scheme with identity section $\xi$, for which $\mu$ becomes a group action extending the regular action.
\end{lemma}
\begin{proof}
    This is clear for $t\neq 0$; for $t=0$ this is \cite[Proposition 2.6]{bai25}.
\end{proof}
Write $\mathcal G=(\mathcal G_t)_{t\in\Omega}\in\operatorname{Coh}(\mathcal M\times_{|\mathcal H|}\mathcal M)$.
\begin{proposition}[Theorem of the square]
    The identity
    \[(\mu\times\operatorname{id})^\ast\mathcal G\cong\operatorname{pr}_{13}^\ast\mathcal G|_{\mathcal M^\circ\times_{|\mathcal H|}\mathcal M}\otimes\operatorname{pr}_{23}^\ast\mathcal G\]
    holds on $\mathcal M^\circ\times_{|\mathcal H|}\mathcal M\times_{|\mathcal H|}\mathcal M$.
\end{proposition}
\begin{proof}
    The fibers of both sides over $\mathcal M^\circ\times_{|\mathcal H|}\mathcal M$ are isomorphic by \cite[Lemma 6.5]{arinkin} and \cite[Theorem 0.3]{bai25}.
    The seesaw principle \cite[Proposition 4.4]{bai25} then shows that they differ by a line bundle pulled back from $\mathcal M^\circ\times_{|\mathcal H|}\mathcal M$, which has to be trivial since their restrictions to the section $\mathcal M^\circ\times_{|\mathcal H|}\mathcal M\times_{|\mathcal H|}\xi$ are both the trivial line bundle.
\end{proof}

\begin{remark}\label{rem:altproof_tots}
    Alternatively, the proposition can be proved by observing that the identity holds on the open set obtained by removing the preimage of the discriminant locus on $|H_0|$.
    Since this locus has codimension $2$ in $|\mathcal H|$, the identity follows as both sides are maximal Cohen-Macaulay.
    This provides an alternative proof of \cite[Theorem 0.3]{bai25} in the case of an elliptic K3 surface.
\end{remark}

Similar to the argument in \cite[\S5]{bai25}, a theorem of the square gives rise to a duality between Tate-Shafarevich twists and Brauer twists.

We refer the reader to \cite{calthesis} for the theory of Brauer twists of derived categories.
Let us briefly introduce the type of Tate-Shafarevich twists we will consider.
Suppose $\varpi:M\to B$ is a proper flat fibration whose generic fiber is an abelian variety.
In the analytic topology, consider the exponential map $\exp:\varpi_\ast T_{M/B}\to\operatorname{Aut}_{M/B}$ for the relative automorphism group.
Let $\operatorname{Aut}^0_{M/B}$ be its image.
Note that this is a sheaf of abelian groups.
\begin{example}
    For the fibration $\mathcal M\to|\mathcal H|$, $\operatorname{Aut}^0_{\mathcal M/|\mathcal H|}$ is precisely the sheaf of sections for the group scheme $\mathcal M^\circ\to |\mathcal H|$.
    The same can be said about $M_t\to |H_t|$, in which case it is the group $P^{\circ\circ}$ in the notation of Kim \cite{kim}.
\end{example}
\begin{definition}
    The Tate-Shafarevich group is the cohomology group $\Sha_{M/B}=H^1(B^{\rm an},\operatorname{Aut}^0_{M/B})$ in the analytic topology.
\end{definition}
Each element $\phi\in\Sha_{M/B}$ defines a new fibration $\varpi^\phi:M^\phi\to B$ via the following procedure:
Represent $\phi=(\phi_{ij})_{ij}$ over an open cover $(U_i)_i$ of $B$.
Then we can reglue the open cover of $M$ given by $(\varpi^{-1}U_i)_i$ with the gluing conditions modified by these automorphisms $\phi_{ij}$.
This gives a new complex manifold $M^\phi$ which, since $\phi_{ij}$ are relative automorphisms, admits a fibration $\varpi^\phi:M^\phi\to B$.

In general, $M^\phi$ has no reason to be projective even if the original fibration is projective.
When $\phi$ is torsion, however, this can be guaranteed.
In fact, these two conditions are essentially equivalent if we restrict ourselves to a particular kind of twist:
Suppose $M\to B=\mathbb P^g$ is a Lagrangian fibration with $M$ a projective irreducible holomorphic symplectic variety.
Then a chosen symplectic form $\sigma$ on $M$ induces an isomorphism $\varpi_\ast T_{M/B}\cong\Omega_B$ (see \cite[Theorem 2.1.11]{aba26}), which allows us to identify
\[H^1(\varpi_\ast T_{M/B})=H^1(\Omega_B)=\mathbb C\]
which, by construction, admits a map to $\Sha_{M/B}$.
\begin{definition}
    When $\varpi$ is a Lagrangian fibration, the space $H^1(\varpi_\ast T_{M/B})$ is called the degenerate twistor line, and its image in $\Sha_{M/B}$ is denoted $\Sha^0_{M/B}$.
\end{definition}
\begin{remark}\label{rem:period_TS}
    We therefore have a family of Lagrangian fibrations parameterized by the degenerate twistor line.
    If $\phi$ is the image of $\eta\in H^1(\Omega_B)$, then the period of $M^\phi$ would be $\sigma+\eta$, where $\sigma$ is the period of $M$.
\end{remark}
\begin{theorem}[{\cite[Theorem A]{aba26}, see also \cite[Theorem B]{AR21}}]\label{thm:torsioniffproj}
    For $\phi\in\Sha^0_{M/B}$, $M^\phi$ is projective if and only if $\phi$ is torsion.
\end{theorem}
\begin{definition}
    We write $\Sha^{\rm alg}_{M/B}$ for the torsion subgroup of $\Sha^0_{M/B}$.
\end{definition}

%If $\phi\in\Sha^{\rm alg}_{M/B}$, then the resulting space $M^\phi$ is algebraic.
\begin{proposition}\label{prop:ts_vs_br}
    Suppose $p\hookrightarrow \Omega$ is either a reduced point or $\Omega$.
    Take $\phi\in\Sha_{\mathcal M_p/|\mathcal H_p|}$.
    Then there is an analytic Brauer class $\alpha_\phi\in\operatorname{Br}^{\rm an}(\mathcal M_p)$ such that:
    \begin{enumerate}
        \item $\mathcal G_p$ is an $(1\boxtimes\alpha_\phi)$-twisted sheaf on $\mathcal M_p^\phi\times_{|\mathcal H_p|}\mathcal M_p$,
        \item the formation of $\alpha_\phi$ is functorial in the sense that for any $p'\hookrightarrow p$, we have $\alpha_{\phi|_{p'}}=\alpha_\phi|_{p'}$, and
        \item $\alpha_{\phi+\phi'}=\alpha_\phi\alpha_{\phi'}$.
    \end{enumerate}
    In particular, if $p=\{t\}\hookrightarrow \Omega$ and $\phi\in\Sha^{\rm alg}_{M_t/|H_t|}$, then $\alpha_\phi$ is torsion, hence is a member of the algebraic Brauer group $\operatorname{Br}(M_t)$.
    In this case $\mathcal G_t$ defines an equivalence
    \[D^b\operatorname{Coh}(M_t^\phi)\cong D^b\operatorname{Coh}(M_t,\alpha_\phi).\]
\end{proposition}

\begin{proof}
    This follows the exact same recipe as \cite[\S5]{bai25}.
    We include a sketch here for completeness.

    The Brauer class $\alpha_\phi$ is constructed as follows:
    Represent $\phi=(\phi_{ij})_{ij}$ over an open cover $(U_i)_i$ of $|\mathcal H_p|$.
    Write $M_{ij}=\mathcal M_p|_{U_{ij}}$ and $M_{ij}^\circ=\mathcal M_p^\circ|_{U_{ij}}$.
    Since $\operatorname{Aut}^0_{\mathcal M_p/|\mathcal H_p|}$ is represented by the group scheme $\mathcal M^\circ_p$, $\phi_{ij}$ may be viewed as a section of $M^\circ_{ij}\to U_{ij}$.
    The observation then is that
    \[(\phi_{ij}\times\operatorname{id})^\ast\mathcal G_p=(\mu\times\operatorname{id})^\ast\mathcal G_p|_{\phi_{ij}\times_{U_{ij}}M_{ij}\times_{U_{ij}}M_{ij}},\]
    which by theorem of the square is nothing but
    \[\operatorname{pr}_2^\ast\mathcal G_p|_{\phi_{ij}\times_{U_{ij}}M_{ij}}\otimes\mathcal G_p.\]

    As $M^\circ_{ij}\to U_{ij}$ is smooth, $\mathcal G_p|_{\phi_{ij}\times_{U_{ij}}M_{ij}}$ is a line bundle.
    The Brauer class $\alpha_\phi$ can then be taken to be the class representing the gerbe $(\mathcal G_p|_{\phi_{ij}\times_{U_{ij}}M_{ij}})_{ij}$.

    The first two properties immediately follow.
    The third property is a consequence of the isomorphism
    \[\mathcal G_p|_{(\phi_{ij}+\phi_{ij}')\times_{U_{ij}}M_{ij}}=\mathcal G_p|_{\phi_{ij}\times_{U_{ij}}M_{ij}}\otimes \mathcal G_p|_{\phi_{ij}'\times_{U_{ij}}M_{ij}}\]
    provided by pulling back theorem of the square to $\phi_{ij}\times_{U_{ij}}\phi_{ij}'\times_{U_{ij}} M_{ij}$.

    Note, finally, that $\alpha_\phi$ is trivial along the fibers of $M_t\to |H_t|$, so the last claim follows from \cite[Theorem 3.2.1]{calthesis}.
\end{proof}
\begin{remark}\label{rem:dualok}
    The same claims would follow if we replace $\mathcal G$ with $\mathcal G^\vee$, simply by dualizing both sides of the identity in theorem of the square.
\end{remark}
\subsection{Constructions of projectively hyperholomorphic kernels}
We now turn to the study of constructions of projectively hyperholomorphic kernels between hyper-K\"ahler manifolds of K3$^{[g]}$-type.
These kernels are possibly twisted vector bundles on the products of hyper-K\"ahler manifolds of K3$^{[g]}$-type, up to multiplying by a line bundle.
Let us refrain from giving a precise definition, as it will not be relevant to our discussion, and refer the reader to \cite{markmanhyperhol,hartliebshah}.

Consider a general polarized K3 surface $(S,H)$ with $\operatorname{Pic}(S)=\mathbb ZH$, and let $n,s$ be integers with $n\ge 2$ and $\gcd(n,s)=1$ so that the Mukai vector
\[v=(n,H,s)\]
is isotropic.
Let $S'=M_S(v)$.
Then $S'$ is another K3 surface and a universal sheaf for this family is a vector bundle $\mathcal U$ on $S\times S'$.
\begin{definition}
    The kernel $(\mathcal U[1])^{[g]+}[-g]=(\mathcal U[-1])^{[g]+}[g]=\mathcal U^{[g]+}\circ\chi^{[g]}=\mathcal U^{[g]}\circ(-\otimes\mathcal O(\delta))$ is called Markman's bundle.
\end{definition}
Markman showed in \cite{markmanhyperhol} that $(\mathcal U[1])^{[g]+}[-g]$ is a projectively hyperholomorphic vector bundle.
\begin{remark}
    In \cite{markmanhyperhol}, the bundle $\mathcal U^{[g]+}$ was instead stated as being projectively hyperholomorphic.
    This seems to be in error due to a sign discrepancy in \cite[Corollary 7.3]{markmanhyperhol}.
    The proof goes through for $(\mathcal U[1])^{[g]+}[-g]$ with a small modification of \cite[Lemma 8.3]{markmanhyperhol}.
\end{remark}

We want to deform Markman's bundle towards a different construction of a projectively hyperholomorphic bundle due to Hartlieb and Shah.
They considered the following situation (\cite[(5.3)]{hartliebshah}):
\[\begin{tikzcd}
    X\drar[swap]{\pi_X}&&Y\dlar{\pi_Y^X}\drar[swap]{\pi_Y^Z}&&Z\dlar{\pi_Z}\\
    &\mathbb P^g&&\mathbb P^g&
\end{tikzcd},\]
where $\pi_X,\pi_Z$ are degree-zero Beauville-Mukai systems of a very general polarized K3 surface of genus $g$, and $\pi_Y^X,\pi_Y^Z$ are twists in $\Sha^{\rm alg}_{X/\mathbb P^g},\Sha^{\rm alg}_{Z/\mathbb P^g}$, respectively, that are transverse Lagrangian fibrations.

Then, correspondingly we have twisted Arinkin transforms giving equivalences
\[\Phi_{\mathcal G_{\rm Arinkin}^\vee}^X:D^b\operatorname{Coh}(Y)\cong D^b(X,\alpha_X),\Phi_{\mathcal G_{\rm Arinkin}^\vee}^Z:D^b\operatorname{Coh}(Y)\cong D^b(Z,\alpha_Z)\]
by \Cref{prop:ts_vs_br}.

They went on to show that the kernel of the composition
\[\Phi_{\mathcal G_{\rm Arinkin}^\vee}^Z\circ(\Phi_{\mathcal G_{\rm Arinkin}^\vee}^X)^{-1}[-g]\]
is a projectively hyperholomorphic (twisted) vector bundle.
\begin{remark}
    The construction we stated is slightly different from the literal version in \cite{hartliebshah}.
    Apart from the obvious manipulations, the dual of Arinkin's sheaf is occasionally taken.
    Let us note that this does not affect the argument in \cite{hartliebshah}.
    We also have $\mathcal G_{\rm Arinkin}^\vee=(\operatorname{id}\times\nu)^\ast\mathcal G_{\rm Arinkin}=(\nu\times\operatorname{id})^\ast\mathcal G_{\rm Arinkin}$ where $\nu$ is the involution given by dualizing a torsion-free sheaf on an integral curve.
\end{remark}
\subsection{Deformation via local Torelli}
We aim to obtain a relationship between these two constructions via a specialization of both towards the same vector bundle on the Hilbert scheme of $g$ points on a K3 surface with two elliptic fibrations, both with integral fibers and a section.

Let $S_0$ be such a surface and $\pi_i,s_i$ the elliptic fibrations and the sections.
We will construct an explicit example later in \Cref{lem:goodk3exists}.

To obtain the desired deformation of spaces, we make use of the local Torelli theorem.
As in \Cref{rem:yoshiokauniv}, write
\[\mathcal D^{[g]}=\{\sigma\in\mathbb P(H^2(S_0^{[g]};\mathbb C)):\sigma^2=0,\sigma\bar\sigma>0\}\]
for the period domain.

For a condition $P(\sigma)$ on the period $\sigma$, we write $\mathcal D^{[g]}|_{P(\sigma)}$ for the locus $\{\sigma\in\mathcal D^{[g]}:P(\sigma)\}$.
Note, then, that for a small period deformation inside $\mathcal D^{[g]}|_{\sigma\cdot f_1=\sigma\cdot f_2=0}$, the fibrations $\pi_i^{[g]}$ both deform \cite{matdeformfib}.
We aim to identify at least one such deformation so that these fibrations on the general fibers are on the degenerate twistor lines of the respective Yoshioka deformations.

Write $\sigma_0$ for the period of $S_0^{[g]}$.
Consider now a Yoshioka deformation of $\pi_1^{[g]}$ and let $\sigma_t^{(1)}$ be the associated period deformation.
In general $f_2$ has no chance of being orthogonal to $\sigma_t^{(1)}$.
But by \Cref{rem:period_TS} we can modify it by a copy of $f_1$ at the cost of performing a Tate-Shafarevich twist.
That is, we are looking for $a_1(t)\in\mathbb C$ so that
\[(\sigma_t^{(1)}+a_1(t)f_1)\cdot f_2=0.\]

Write $n=f_1\cdot f_2$.
Then this solves to
\[a_1(t)=-\frac{1}{n}\sigma_t^{(1)}\cdot f_2,\]
and the period deformation $\sigma_t=\sigma_t^{(1)}+a_1(t)f_1$ locally preserves both fibrations.

We aim to realize $\sigma_t$ as a similar Tate-Shafarevich twist of a Yoshioka deformation of $\pi_2^{[g]}$, at least some of the times.
This requires some $a_2(t)$ so that
\[(\sigma_t-a_2(t)f_2)\cdot s_2^+=0,s_2^+=\left(s_2-\frac{\delta}{2}\right),\]
which solves to
\[a_2(t)=\sigma_t\cdot s_2^+.\]

Let us write $\sigma_t^{(2)}=\sigma_t-a_2(t)f_2$.
Note that both $a_i(t)$ are holomorphic functions with $a_i(0)=0$.
\begin{lemma}\label{lem:reconstr_hartliebshah}
    For a general choice of the initial deformation $\sigma_t^{(1)}$, the deformation $\sigma_t^{(2)}$ is a Yoshioka deformation.
\end{lemma}
\begin{proof}
    Essentially, $\sigma_t^{(2)}$ can be obtained from $\sigma_t^{(1)}$ via the linear map given by
    \[\sigma\mapsto \sigma-\frac{1}{n}(\sigma\cdot f_2)f_1-\left(\sigma\cdot s_2^+-\frac{1}{n}(\sigma\cdot f_2)(f_1\cdot s_2^+)\right)f_2\]
    by going through the steps in our procedure.
    This map gives an isomorphism
    \[\mathcal D^{[g]}|_{\sigma\cdot f_1=\sigma\cdot (s_1-\delta/2)=0}\to \mathcal D^{[g]}|_{\sigma\cdot f_2=\sigma\cdot (s_2-\delta/2)=0}.\]

    By \Cref{rem:yoshiokauniv}, a very general local deformation in $\mathcal D^{[g]}|_{\sigma\cdot f_i=\sigma\cdot (s_i-\delta/2)=0}$ at $\sigma_0$ is an instance of a Yoshioka deformation for $\pi_i^{[g]}$.
\end{proof}
Along the period deformation $\sigma_t$, we then have a family over a small disk $\Omega\ni t$:
\[\begin{tikzcd}
    M_t^{(1)}\drar&&Y_t\drar\dlar&&M_t^{(2)}\dlar\\
    &\mathbb P^g&&\mathbb P^g&
\end{tikzcd}\]
where $M_t^{(i)}\to \mathbb P^g$ are Yoshioka deformations of $\pi_i^{[g]}$, and $Y_t$ corresponds to the period deformation given by $\sigma_t$, which by construction admits Lagrangian fibrations corresponding to the classes $f_i$ that are Tate-Shafarevich twists in $\Sha^{\rm alg}_{M_t^{(i)}/\mathbb P^g}$.
\begin{lemma}
    Suppose $Y$ is hyper-K\"ahler with two Lagrangian fibrations $Y\to\mathbb P^g$ with classes $f_1,f_2$ with $f_1\cdot f_2>0$.
    Then $Y$ is projective.
    If $f_1$ has integral fibers, then the two fibrations are transverse in the sense that the combined map $Y\to\mathbb P^g\times\mathbb P^g$ is finite.
\end{lemma}
The condition $f_1\cdot f_2>0$ is just saying that $f_1,f_2$ represent distinct fibrations.
\begin{proof}
    Let us show that $f_2$ is relatively ample for $f_1$.
    This suffices since if there were a pair of fibers $F_1,F_2$ for $f_1,f_2$ intersecting at positive dimensions, then one can take a curve $C$ inside their intersection.
    This curve has the property that $C\cdot f_2=0$ since it is contained in $F_2$.
    But it is also contained in $F_1$, so $f_2$ cannot be ample on $F_1$.

    The class $f_1+f_2$ is positive, so $Y$ is projective \cite[Theorem 3.11]{huybasic}.
    Let $H$ be an ample class on $Y$.
    If $F_1$ is a smooth fiber of $f_1$, then
    \[\lambda H|_{F_1}=f_2|_{F_1}\]
    for some $\lambda\in\mathbb Q$ by \cite{matdeformfib}.

    By \cite[Lemma 4.4]{voisintwistedlsv}, $\lambda H|_{F_1}=f_2|_{F_1}$ numerically for every fiber $F_1$ of $f_1$.
    The condition $f_1\cdot f_2>0$ implies $\lambda>0$, hence $f_2$ is relatively ample for $f_1$.
\end{proof}
We have therefore reconstructed the situation \cite[(5.3)]{hartliebshah} without having to make use of the SYZ conjecture for K3$^{[g]}$-type.
This construction now allows us to take a limit of Hartlieb-Shah's kernel.

Write $\mathcal M^{(i)}=(M_t^{(i)})_{t\in\Omega}$ and $\mathcal Y=(Y_t)_{t\in\Omega}$.
Denote by $\mathcal G^{(i)}$ the deformation of kernel as considered before for the Yoshioka deformation $\mathcal M^{(i)}\to\mathbb P^g_\Omega$.
Let us now apply \Cref{prop:ts_vs_br}.
\begin{proposition}
    There is a Brauer class $\alpha^{(i)}\in\operatorname{Br}^{\rm an}(\mathcal M^{(i)})$ such that:
    \begin{enumerate}
        \item $\mathcal G^{(i)}$ is a $(1\boxtimes\alpha^{(i)})$-twisted sheaf on $\mathcal Y\times_\Omega\mathcal M^{(i)}$;
        \item for each $t\neq 0$, $\alpha^{(i)}_t=\alpha^{(i)}|_{M_t^{(i)}}$ is algebraic, i.e.~it lives inside $\operatorname{Br}(M_t^{(i)})$, and $\mathcal G^{(i)}_t$ is the twisted Arinkin transform
        \[\Phi_{\mathcal G_{\rm Arinkin}^\vee}:D^b\operatorname{Coh}(Y_t)\to D^b\operatorname{Coh}(M_t^{(i)},\alpha^{(i)}_t);\]
        as considered by Hartlieb and Shah;
        \item $\alpha^{(i)}_0=0$, and $\mathcal G^{(i)}_0$ is the BKR congruence $\mathcal P_i^{[g]}$, where $\mathcal P_i$ is the kernel for the elliptic fibration $\pi_i$.
    \end{enumerate}
\end{proposition}
\begin{proof}
    Since $a_i(t)$ is holomorphic and $a_i(0)=0$, we have a family of degenerate twistor deformations.
    More precisely, writing $\varpi:\mathcal M^{(i)}\to\mathbb P^g_\Omega$, we have a global section of
    \[\mathbf R^1 b_\ast(\varpi_\ast T_{\mathcal M^{(i)}/\mathbb P^g_\Omega})=\mathbf R^1 b_\ast\Omega_{\mathbb P^g_\Omega/\Omega}\]
    where $b:\mathbb P^g_\Omega\to\Omega$.
    As $\Omega$ is a small disk, this gives us an element of
    \[H^1(\varpi_\ast T_{\mathcal M^{(i)}/\mathbb P^g_\Omega})\cong H^0(\mathbf R^1 b_\ast(\varpi_\ast T_{\mathcal M^{(i)}/\mathbb P^g_\Omega}))\]
    and thus, via the exponential map, an element in $\Sha_{\mathcal M^{(i)}/\mathbb P^g_\Omega}$ restricting to a degenerate twistor line deformation at each fiber.
    \Cref{prop:ts_vs_br} then yields the result.
    The algebraicity of the Brauer classes follows from the fact that $Y_t$ is projective due to \Cref{thm:torsioniffproj}.
\end{proof}
The kernel constructed by Hartlieb-Shah is then the generic fiber of
\[\mathcal G^{(1)}\circ(\mathcal G^{(2)})^{-1}[-g],\]
and hence:
\begin{corollary}\label{cor:hartliebshah_specialize}
    The kernel $\mathcal P_1^{[g]}\circ(\mathcal P_2^{[g]})^{-1}[-g]$ is a specialization of Hartlieb-Shah's construction.
\end{corollary}
We will show that an instance of Markman's bundle can be deformed to where it can be factored into the form in the Corollary.
\subsection{A K3 surface with two elliptic fibrations}
Let us finally construct a K3 surface $S_0$ with two elliptic fibrations $\pi_1,\pi_2:S_0\to\mathbb P^1$, both with integral fibers and a section.
\begin{example}
    Any line on a smooth quartic defines an elliptic fibration by projecting from it.
    Take a quartic with only two lines that are skew, then the two lines give two elliptic fibrations, both of which have to have integral fibers due to the lack of other lines.
    Both also have a section, since one line would provide a section for the fibration induced by the other.
\end{example}
In this example, the intersection number between general fibers of the two fibrations is $2$.
To construct one where the intersection number is any $n\ge 2$, let us consider the lattice $\mathbb Zf_1\oplus\mathbb Zf_2\oplus\mathbb Zs$ under the pairing
\[\begin{pmatrix}
    0&n&1\\
    n&0&1\\
    1&1&-2
\end{pmatrix}.\]
Note that this quadratic form has signature $(1,2)$.
\begin{lemma}\label{lem:goodk3exists}
    There is a projective K3 surface $S_0$ whose N\'eron-Severi lattice is exactly the lattice $\mathbb Zf_1\oplus\mathbb Zf_2\oplus\mathbb Zs$ described above.
    Moreover, $S_0$ can be chosen so that
    \begin{enumerate}
        \item $f_1+f_2$ is ample,
        \item $f_1,f_2,s$ are all represented by smooth curves, and
        \item the elliptic fibrations $\pi_i$ induced by $f_i$ have integral fibers.
    \end{enumerate}
\end{lemma}
In particular, $s$ would be a section to both the elliptic fibrations $\pi_i:S_0\to\mathbb P^1$.
\begin{proof}
    By \cite[Ch.~14, Corollary 3.1]{huyk3}, there is a K3 surface with this N\'eron-Severi lattice.
    It remains to show that this K3 surface can be chosen so that the assertions are true.

    Let us first show that $h=f_1+f_2$ can be made ample.
    By \cite[Ch.~8, Corollary 2.9]{huyk3}, it suffices to show that it has nonzero intersection with any $(-2)$-class.

    So let $c=xf_1+yf_2+zs$ be such that $c^2=-2$, with $x,y,z\in\mathbb Z$.
    Then $nxy+xz+yz-z^2=-1$, in particular $z\neq 0$.
    Now suppose $c\cdot h=0$, then $c\cdot f_1=-c\cdot f_2$.
    On the other hand,
    \begin{align*}
        (c\cdot f_1)(c\cdot f_2)&=(ny+z)(nx+z)=n^2xy+nxz+nyz+z^2\\
        &=n(z^2-1)+z^2>0,
    \end{align*}
    so this is impossible.
    Therefore $h$ is ample.

    Now suppose $c=xf_1+yf_2+zs$ is the class of an irreducible curve.
    We want to show that $c\cdot f_1\ge0,c\cdot f_2\ge 0$.
    Note that $c\cdot f_1+c\cdot f_2=c\cdot h>0$.
    On the other hand,
    \[(c\cdot f_1)(c\cdot f_2)=n\left(z^2+\frac{c^2}{2}\right)+z^2=(n+1)z^2+n\frac{c^2}{2}.\]
    If $c^2=-2$, then as before we have $z\neq 0$ and therefore this product is positive; if $c^2\neq -2$, then $c^2\ge 0$, so this product is nonnegative.
    Hence $c\cdot f_1\ge0,c\cdot f_2\ge 0$.
    So $f_1,f_2$ are nef.

    As for $s$, note that it is effective since $s\cdot h=2>0$.
    Suppose
    \[s=\left[\sum_{i=1}^Nm_iC_i\right]\]
    for $C_i$ irreducible and $m_i>0$.
    Since $s^2<0$, we may assume WLOG that $C_1$ has negative square.
    Then $C_1^2=-2$ and our previous calculation shows $C_1\cdot f_j>0$ for $j=1,2$.

    But we also have $s\cdot f_j=1$.
    Applying this to $f_1$ gives $m_1=1$ and $[C_i]=f_1$ for any $i>1$.
    Then doing the same to $f_2$ gives $[C_i]=f_2$ for any $i>1$, which is absurd unless $N=1$.
    Hence any effective representative of $s$ is already reduced and irreducible.
\end{proof}

Write $\mathcal P_i$ for the autoequivalence kernel for the elliptic K3 fibration $\pi_i$ as considered before.
\begin{lemma}
    The composition $\mathcal E=\mathcal P_1\circ\mathcal P_2^\vee$ is a vector bundle of rank $n$, and there is a polarization for which, for every $x\in S_0$, $\Phi_{\mathcal E}(k(x))$ is slope-stable.
\end{lemma}
\begin{proof}
    This is essentially \cite[Theorem 6.6]{bottini}.
\end{proof}
\begin{lemma}
    There exists a line bundle $L\in\operatorname{Pic}(S_0)$ such that $\mathcal E\otimes\operatorname{pr}_2^\ast L$ is a universal bundle for $S_0\cong M_{S_0}(n,L_0,n-1)$ where $L_0$ is a primitive ample class.
\end{lemma}
\begin{proof}
    By the preceding lemma, it suffices to check that there is some $L$ such that $v(\Phi_{\mathcal E}(k(x))\otimes L)=(n,H,n-1)$ for a primitive ample $H$.

    Let us first compute $v(\Phi_{\mathcal E}(k(x)))$.
    Note that $v(\Phi_{\mathcal P_2^\vee}(k(x)))=f_2$.
    By \Cref{prop:elliptic_surface_calc}, we have
    \begin{align*}
        v(\Phi_{\mathcal E}(k(x)))&=\Phi_{\mathcal P_1}^{\tilde{H}}(f_2)\\
        &=\Phi_{\mathcal P_1}^{\tilde{H}}(f_2-ns-(2n+1)f_1)+n\Phi_{\mathcal P_1}^{\tilde{H}}(s)+(2n+1)\Phi_{\mathcal P_1}^{\tilde{H}}(f_1)\\
        &=-f_2+ns+(2n+1)f_1+n(\alpha-f_1+\beta)-(2n+1)\beta\\
        &=n\alpha+(ns+(n+1)f_1-f_2)+(-n-1)\beta.
    \end{align*}

    So we choose $L=\mathcal O(-s-f_1+f_2)$.
    Then
    \[v(\Phi_{\mathcal E}(k(x))\otimes L)=n\alpha+(f_1+(n-1)f_2)+(n-1)\beta,\]
    and $f_1+(n-1)f_2=f_1+f_2+(n-2)f_2$ is primitive and ample.
\end{proof}
By deforming the pair $(S_0,L_0)$, we immediately obtain from \cite[Proposition 5.1]{yoshioka} that:
\begin{corollary}
    There are line bundles $L,L'$ such that $\mathcal E'=\mathcal E\otimes (L\boxtimes L')$ is the central fiber of a family of universal bundles $\mathcal U_t$ on $S_t\times M_{S_t}(n,L_t,n-1)$ where $\operatorname{Pic}S_\eta=\mathbb ZL_\eta$.
\end{corollary}

\begin{proposition}\label{prop:markmanspecialize}
    The object $\mathcal P_1^{[g]}\circ(\mathcal P_2^{[g]})^{-1}[-g]$ is a rigid, projectively hyperholomorphic bundle and, up to a line bundle, is a specialization of Markman's bundle.
\end{proposition}
\begin{proof}
    Observe that
    \begin{align*}
        \mathcal P_1^{[g]}\circ(\mathcal P_2^{[g]})^{-1}&=\operatorname{BKR}^\top\circ\mathcal P_1^{\boxtimes g}\circ\operatorname{BKR}\circ\operatorname{BKR}^{-1}\circ(\mathcal P_2^{-1})^{\boxtimes g}\circ(\operatorname{BKR}^\top)^{-1}\\
        &=(\mathcal P_1\circ\mathcal P_2^{-1})^{[g]+}\\
        &=((-\otimes L_2)\circ\mathcal E'\circ (-\otimes L_1)[1])^{[g]+}\\
        &=(-\otimes L_2)^{[g]+}\circ(\mathcal E'[1])^{[g]+}\circ (-\otimes L_1)^{[g]+}
    \end{align*}
    for some line bundles $L_1,L_2$ on $S_0$.

    But note that $(-\otimes L)^{[g]+}$ is again tensoring by a line bundle.
    Indeed, $\operatorname{BKR}^\top(A\otimes L^{\boxtimes g})=\operatorname{BKR}^\top(A)\otimes \operatorname{HC}^\ast L^{(g)}$ by \cite[Remark 3.11]{krug}, where $L^{(g)}$ is the symmetric power of $L$ on $S_0^{(g)}$, and $\operatorname{HC}:S_0^{[g]}\to S_0^{(g)}$ is the Hilbert-Chow resolution.

    So, up to a line bundle, the composition is $(\mathcal E'[1])^{[g]+}[-g]=(\mathcal E')^{[g]}\circ(-\otimes\mathcal O(\delta))$, which is a projectively hyperholomorphic bundle by the argument of \cite{markmanhyperhol}.
    It is also rigid by \Cref{lem:FMkernelsarerigid}.
    By the preceding corollary, it is a specialization of Markman's bundle.
\end{proof}

\bibliographystyle{plain}
\bibliography{ref.bib}
\end{document}